\pdfoutput=1

\documentclass[final,3p,times]{elsarticle}

\usepackage[utf8]{inputenc}
\usepackage{bm}
\usepackage{graphicx}
\usepackage{amsmath}
\usepackage{amssymb}
\usepackage{amstext}
\usepackage{amsfonts}
\usepackage[labelsep=period,textfont={bf}]{caption}
\usepackage[labelformat=simple]{subcaption}
\usepackage{epstopdf}
\usepackage{comment}
\usepackage{soul}
\usepackage[colorlinks=true, pdfstartview=FitV, linkcolor=black, citecolor= black, urlcolor= black]{hyperref}
\usepackage{footnote}
\usepackage{footnpag}			      	
\usepackage{multirow}
\usepackage[table,xcdraw]{xcolor}
\usepackage{siunitx}
\usepackage{lscape}
\usepackage[nolist,nohyperlinks]{acronym}
\usepackage{booktabs}
\usepackage{tikz}
\usetikzlibrary{automata,positioning}
\usetikzlibrary{shapes.multipart}
\usetikzlibrary{shapes.geometric}
\setcitestyle{square,comma}

\begin{acronym}
\acro{CR3BP}{Circular Restricted Three Body Problem}
\acro{TPBVP}{Two-Point Boundary-Value Problem}
\acro{PSO}{Particle Swarm Optimization}
\acro{STM}{State Transition Matrix}
\acro{PMP}{Pontryagin's Maximum Principle}
\acro{COV}{Calculus of Variations}
\acro{ACT}{Adjoint Control Transformation}
\end{acronym}

\journal{Acta Astronautica}

\begin{document}

\begin{frontmatter}

\renewcommand\thesubfigure{\alph{subfigure}} 

\title{Particle Swarm Optimization-Based Co-State Initialization for Low-Thrust Minimum-Fuel Trajectory Optimization}

\author[label1]{Grant R. Hecht}
\ead{granthec@buffalo.edu}

\author[label1]{Eleonora M. Botta\corref{corl}}
\cortext[corl]{Corresponding author.}
\ead{ebotta@buffalo.edu}

\address[label1]{Department of Mechanical and Aerospace Engineering, University at Buffalo, Buffalo, NY, United States}

\begin{abstract}
In this paper, Particle Swarm Optimization with energy-to-fuel continuation is proposed for initializing the co-state variables for low-thrust minimum-fuel trajectory optimization problems in the circular restricted three-body problem. Particle Swarm Optimization performs a search of the solution space by minimizing the weighted sum of squares of the two-point boundary-value problem final boundary condition residuals for the minimum-energy problem. Next, an energy-to-fuel homotopy is employed to transition the minimum-energy trajectory to a minimum-fuel trajectory, starting from the generated guess. The proposed methodology is applied to two low-thrust transfer problems in the Earth-Moon system: a transfer from a geostationary transfer orbit to an L1 halo orbit, and a transfer from an L2 halo orbit to an L1 halo orbit. The resulting minimum-fuel trajectories are validated with the literature. It is demonstrated that the methodology can successfully generate guesses for the initial co-state variables which converge to a solution for both scenarios. A strategically chosen particle swarm size is shown to improve the convergence rate of the methodology. The proposed approach is of simple implementation, can be easily extended to other trajectory optimization problems, facilitates the discovery of multiple candidate trajectories, and does not require a user-provided guess, all of which are advantageous features for the preliminary phase of mission design.
\end{abstract}


\end{frontmatter}

\section{Introduction}
Low-thrust propulsion has gained much attention in recent years and is now often considered the best option for a wide range of space missions. The associated reduction in required propellant mass compared to traditional chemical propulsion results in cost savings and an increase in the feasible payload ratio. Low-thrust propulsion does, however, require long, continuous thrusting arcs and trajectories often spanning many revolutions about a large body, which can lead to high-dimensional optimization problems with many local minima. Low-thrust spacecraft trajectory optimization problems most commonly focus on minimizing the time of flight, fuel required, or combinations of these conflicting objectives. As the name suggests, minimum-time trajectory optimization is necessary when mission requirements specify that the spacecraft must reach its destination in a timely manner, regardless of the fuel required, as is often the case when transporting human life or trying to acquire time-sensitive scientific data. On the other hand, minimum-fuel trajectory optimization is useful when time is not such an important constraint and it is desired to reduce the required fuel, which allows for increasing the feasible payload ratio and decreasing mission costs. Many approaches for optimizing low-thrust spacecraft trajectories have been proposed and applied successfully. Typically, approaches can be separated into two categories known as \textit{direct} and \textit{indirect} methods. Hybrid approaches also exist which combine the characteristics of direct and indirect methods \cite{Woollands_2019}. A broad review of the state-of-the-art techniques employed for spacecraft trajectory optimization can be found in work by Chai et al. \cite{CHAI_2019}.

Direct trajectory optimization methods formulate the trajectory optimization problem by parameterizing the state and control variables, forming a nonlinear programming (NLP) problem \cite{Woollands_2019}. These methods can be further categorized into shooting, collocation, dynamic programming, or differential inclusion-based methods \cite{CHAI_2019}. When considering trajectories spanning tens to hundreds of days with many separate, continuous thrusting arcs -- as is common for low-thrust propulsion --, this parameterization can lead to very high-dimensional problems which may be intractable if not strategically formulated, due to the approaching curse of dimensionality. Additionally, this parameterization generally enforces the control to take some predetermined form, which can result in suboptimal solutions. Nevertheless, direct methods continue to be an attractive approach for solving low-thrust trajectory optimization problems for several reasons. When compared to indirect methods, direct methods can more easily incorporate complicated constraints and tend to be much less sensitive to the quality of the initial guess supplied to the algorithm \cite{Woollands_2019}. Furthermore, physical intuition can be exploited when an initial guess for the solution is determined and supplied to an NLP solver.

Indirect methods apply \ac{COV} and \ac{PMP} to analytically derive the necessary conditions for optimality, introducing additional co-state variables and formulating the problem as a \ac{TPBVP} \cite{Taheri_2017}. Collocation or finite element-based methods can then be used to approximate the solution of the \ac{TPBVP}, or a shooting-based method may be used to solve the \ac{TPBVP} directly \cite{CHAI_2019}. Collocation or finite element-based approaches can lead to useful physical insight about the problem while providing approximate solutions \cite{Woollands_2019}. Alternatively, shooting-based methods result in a fully continuous optimal control law for the low-thrust spacecraft trajectory. Unfortunately, the convergence of the indirect shooting-based approach is typically not guaranteed and is highly dependent on the selection of the initial co-state variables used to initialize the algorithm \cite{Woollands_2019,Taheri_2017}.

Typically, the co-state variables are abstract quantities with no clear physical significance and the convergence radius of indirect shooting methods is often very small. Therefore, determining an initial choice of co-state variables to initialize a shooting routine can be a daunting task. To mitigate some of this difficulty, several methods have been proposed to assist in initializing the co-state variables. One such approach, first proposed by Dixon and Biggs \cite{Dixon_1972,Dixon_1981} and termed \ac{ACT}, involves relating the co-state variables to alternative variables of physical significance, thereby decreasing the sensitivity of the problem and allowing for the use of physical insight when selecting a guess. Application of \ac{ACT} to problems in astrodynamics is discussed, for example, in References \cite{Senent_2005,Ranieri_2005,Ranieri_2006,Russell_2007}. An alternative approach involves formulating an approximation of the initial co-state variables to avoid the need for a user-provided guess. For example, within the field of astrodynamics, Thorne and Hall approximated the initial co-state variables by solving a similar problem with an analytical solution by assuming zero gravity and a constant spacecraft mass \cite{Thorne_1996}, whereas Lee et al. derived an approximate relationship between the co-state variables and the initial constrained physical state based on the behavior of the initial co-states for a range of optimal trajectories \cite{Lee_2009,Lee_2012}. Unfortunately, these techniques are restricted to planar motion or transfers of less than one revolution. Another possibility to estimate the initial co-state variables involves employing a solution acquired via direct trajectory optimization. Examples of such methods include a technique that exploits approximate solutions acquired via shape-based methods and nonlinear least squares \cite{Taheri_2017}, and an approach that allows one to relate the Lagrange multipliers of the Karush-Kuhn-Tucker equations associated with various collocation methods to the co-states of the optimal control problem \cite{Fahroo_2001,Benson_2006,Garg_2011,Francolin_2015,Eide_2018}. Of course, a primary difficulty of such approaches is that a solution must be found through the appropriate direct trajectory optimization method first, a task that can itself be difficult. Furthermore, once a solution is discovered via direct trajectory optimization, successful estimation of the co-state variables is not guaranteed.

A particularly interesting approach, proposed by Pontani and Conway, employs \ac{PSO} to select co-state variables which minimize the sum of the absolute value of the final boundary condition residuals \cite{Pontani_2010,Pontani_2014}, solving the \ac{TPBVP} directly. Such an approach is highly desirable as it can, in principle, be applied to any \ac{TPBVP} and does not require a user to provide any guess (although box constraints must be defined to limit the search space). Unfortunately, due to the extreme sensitivity of the final boundary condition residuals to perturbations in the initial co-state variables for many indirect trajectory optimization problems, this initial application of \ac{PSO} was limited to planar minimum-time transfers in the two-body problem, which reduces the number of unknown variables and constraints. 

In this paper, a method is proposed which builds upon the work by Pontani and Conway \cite{Pontani_2010,Pontani_2014} and allows for solving non-planar low-thrust minimum-fuel trajectory optimization problems with the indirect shooting approach by leveraging \ac{PSO} to initialize the co-state variables. \ac{PSO} co-state initialization and an energy-to-fuel homotopy technique are coupled \cite{Bertrand_2002} to widen the convergence radius of the \ac{TPBVP} and reduce the difficulty of the parameter optimization problem to be solved by \ac{PSO}. \ac{PSO} generates co-states for the easier-to-solve minimum-energy problem, before transitioning the solution to minimum-fuel through continuation, all without the need of a user-provide guess. Herein, the proposed methodology is presented in the framework of the \ac{CR3BP}, and important considerations in regard to solution optimality and computational performance are discussed. Case studies, based on previous works \cite{Zhang_2015,Aziz_2019}, are presented which shed light on the promising performance of the proposed methodology while allowing for its validation. Guidance on the selection of the particle swarm size, a heuristic parameter of the optimization algorithm which has a significant effect on co-state initialization performance, is also provided.  

The remainder of this paper is organized as follows. In Section \ref{sec:formulation}, the problem is formulated, beginning with the \ac{CR3BP} dynamics, and then deriving the indirect minimum-fuel \ac{TPBVP} with an energy-to-fuel homotopy through the application of \ac{COV} and \ac{PMP}. The proposed method of co-state initialization is discussed in Section \ref{sec:TPBVP}, along with the single shooting algorithm used to robustly transition \ac{PSO}-initialized co-state values, corresponding to the minimum-energy problem, to a solution of the minimum-fuel problem. In Section \ref{sec:results}, two minimum-fuel trajectory optimization scenarios (i.e., a transfer from a geostationary transfer orbit to L1 halo orbit, and a transfer from an L2 halo orbit to an L1 halo orbit)
are employed to validate the methodology, demonstrate the performance of the proposed method, and analyze the effect of the particle swarm size on convergence. Finally, conclusions are made which summarize the achievements of this work and future steps. 

\section{Indirect Optimal Control Formulation}
\label{sec:formulation}
In this Section, the indirect optimal control \ac{TPBVP} is formulated for the minimum-fuel problem with terminal state constraints and a fixed time of flight. This is done using the \ac{CR3BP} dynamics model, which is provided and discussed. Throughout the \ac{TPBVP} formulation, an energy-to-fuel homotopy technique is incorporated by introducing a perturbed energy term to the minimum-fuel cost function, before the Hamiltonian is formed and the analytical necessary conditions of optimality are derived.

\subsection{Circular Restricted Three Body Problem Dynamics}
The \ac{CR3BP} dynamics describe the motion of a body in a synodic (rotating) reference frame under the influence of two primary bodies which travel along co-planar circular paths centered about the barycenter of the system \cite{SZEBEHELY_1967_Chap1}. The mass of the third body is indicated with $m$ and is assumed to be much lower than those of the primary bodies $m_1$ and $m_2$; therefore it has no effect on their motion. In the following, all variables are non-dimensionalized such that the angular velocity of both primary bodies, the sum of the masses of the primary bodies, and the distance between the primary bodies are set to unity to improve numerical stability. If the third body is a spacecraft with the capability to apply propulsive thrust in any direction, the equations of motion describing the evolution of the spacecraft's position and velocity in a synodic reference frame, along with its mass are given by \cite{Russell_2007}
\begin{align}
    \Dot{\mathbf{x}} = \mathbf{f}(\mathbf{x},\Hat{\boldsymbol{\alpha}},u) = \begin{bmatrix}
    \Dot{\mathbf{r}} \\ \Dot{\mathbf{v}} \\ \Dot{m}
    \end{bmatrix} = \begin{bmatrix}
    \mathbf{v} \\ \mathbf{g}(\mathbf{r}) + \mathbf{h}(\mathbf{v}) + uT_{max}\hat{\boldsymbol{\alpha}}/m \\ -uT_{max}/c
    \end{bmatrix}
    \label{eqn:state_dyn}
\end{align}

\noindent where $\mathbf{r}$ and $\mathbf{v}$ are the non-dimensionalized position and velocity vectors of the spacecraft, $T_{max}$ is its maximum capable thrust, and $c$ is the propellant exhaust velocity, which is assumed to be constant and given by $c=I_{sp}g_0$, with $I_{sp}$ the specific impulse of the propellant and $g_0=9.81$ m/s$^2$. The control variables are $u$ and $\hat{\boldsymbol{\alpha}}$, where $u \in [0,1]$ is the thrust throttling factor and $\hat{\boldsymbol{\alpha}}$ is a unit vector representing the direction at which thrust is applied. Functions $\mathbf{g}(\mathbf{r})$ and $\mathbf{h}(\mathbf{v})$ are given by \cite{Zhang_2015}
\begin{align}
    \mathbf{g}(\mathbf{r}) &= \begin{bmatrix}
        r_x - \frac{(1 - \mu)(r_x + \mu)}{r_1^3} - \frac{\mu(r_x + \mu - 1)}{r_2^3} \\
        r_y - \frac{(1 - \mu)r_y}{r_1^3} - \frac{\mu r_y}{r_2^3} \\
        -\frac{(1 - \mu)r_z}{r_1^3} - \frac{\mu r_z}{r_2^3}
    \end{bmatrix} \\ 
    \mathbf{h}(\mathbf{v}) &= \begin{bmatrix} 2v_y & -2v_x & 0 \end{bmatrix}^T
\end{align}
\noindent where subscripts $x$, $y$, and $z$ denote the components of the position and velocity vector along the corresponding Cartesian direction of the synodic reference frame and $\mu = m_2/(m_1 + m_2)$ is the non-dimensionalized mass of the second primary body. Additionally, $r_1$ and $r_2$ denote the non-dimensionalized distance of the spacecraft from the first and second primary bodies respectively and are given by
\begin{align}
    r_1 &= \sqrt{(r_x + \mu)^2 + r_y^2 + r_z^2} \\
    r_2 &= \sqrt{(r_x + \mu - 1)^2 + r_y^2 + r_z^2}
\end{align}

\subsection{Indirect Minimum-Fuel Trajectory Optimization using Energy-to-Fuel Homotopy}
For the minimum-fuel problem, it is desired to find a solution of Eq. \eqref{eqn:state_dyn} which minimizes the cost function 
\begin{align}
    J_{mf} = \frac{T_{max}}{c}\int_{t_i}^{t_f}udt
\end{align}
\noindent where $t_i$ and $t_f$ correspond to the initial and final time of the trajectory. A known characteristic of minimum-fuel trajectories is the bang-bang nature of their control profile \cite{Russell_2007,Bertrand_2002,Zhang_2015}, which results in discontinuities in the throttling factor $u$ along the trajectory and can severely restrict the convergence radius of numerical methods. To mitigate this issue, Bertrand and Epenoy proposed the use of an \textit{energy-to-fuel} homotopy continuation for transforming solutions from the discontinuity-free minimum-energy problem to the minimum-fuel problem \cite{Bertrand_2002}. A perturbed energy term is introduced into the cost function by means of the perturbation parameter $\epsilon$, which is gradually reduced until its effect is nulled. Employing this approach, the new cost function is given by
\begin{align}
    J = \frac{T_{max}}{c}\int_{t_i}^{t_f}\left[u - \epsilon u(1 - u)\right]dt \hspace{5mm} \epsilon \in [0, 1]
    \label{eqn:homotopycost}
\end{align}

To begin the derivation of the \ac{TPBVP}, we proceed by forming the Hamiltonian as \cite{Zhang_2015}
\begin{align}
    H &= \boldsymbol{\lambda}^T\mathbf{f}(\mathbf{x},\Hat{\boldsymbol{\alpha}},u) + \frac{T_{max}}{c}\left[u - \epsilon u(1 - u)\right]\\
      &= \boldsymbol{\lambda}_r^T\mathbf{v} + \boldsymbol{\lambda}_v^T\left[\mathbf{g}(\mathbf{r}) + \mathbf{h}(\mathbf{v}) + \frac{uT_{max}}{m}\Hat{\boldsymbol{\alpha}}\right] - \lambda_m\frac{uT_{max}}{c} + \frac{T_{max}}{c}\left[u - \epsilon u(1 - u)\right]
      \label{eqn:H1}
\end{align}
\noindent where $\boldsymbol{\lambda}^T = \left[\boldsymbol{\lambda}_r^T,\boldsymbol{\lambda}_v^T,\lambda_m \right]$ is the introduced co-state vector, which evolves in time according to
\begin{align}
    \Dot{\boldsymbol{\lambda}} = \begin{bmatrix} \Dot{\boldsymbol{\lambda}}_r \\ \Dot{\boldsymbol{\lambda}}_v \\ \Dot{\lambda}_m
    \end{bmatrix} = -\left(\frac{\partial H}{\partial \mathbf{x}}\right)^T = \begin{bmatrix} -\mathbf{G}^T\boldsymbol{\lambda}_v \\ -\boldsymbol{\lambda}_r - \mathbf{H}^T\boldsymbol{\lambda}_v \\ \frac{uT_{max}}{m^2}\boldsymbol{\lambda}_v^T\hat{\boldsymbol{\alpha}}
    \end{bmatrix}
    \label{eqn:costate_dyn}
\end{align}
\begin{align}
    \mathbf{G} &= \frac{\partial \mathbf{g}(\mathbf{r})}{\partial \mathbf{r}} \\
    \mathbf{H} &= \frac{\partial \mathbf{h}(\mathbf{v})}{\partial \mathbf{v}}
\end{align} 
\noindent where the non-zero terms of $\mathbf{G}$ and $\mathbf{H}$ are
\begin{align}
    G_{1,1} &= 1 - \frac{1 - \mu}{r_1^3} + \frac{3(1-\mu)(x + \mu)^2}{r_1^5} - \frac{\mu}{r_2^3} + \frac{3\mu(x + \mu - 1)^2}{r_2^5} \\
    G_{2,2} &= 1 - \frac{1 - \mu}{r_1^3} + \frac{3(1-\mu)y^2}{r_1^5} - \frac{\mu}{r_2^3} + \frac{3\mu y^2}{r_2^5} \\
    G_{3,3} &= -\frac{1 - \mu}{r_1^3} + \frac{3(1-\mu)z^2}{r_1^5} - \frac{\mu}{r_2^3} + \frac{3\mu z^2}{r_2^5} \\
    G_{1,2} &= G_{2,1} = \frac{3(1-\mu)(x+\mu)y}{r_1^5} + \frac{3\mu(x + \mu - 1)y}{r_2^5} \\
    G_{1,3} &= G_{3,1} = \frac{3(1-\mu)(x+\mu)z}{r_1^5} + \frac{3\mu(x + \mu - 1)z}{r_2^5} \\
    G_{2,3} &= G_{3,2} = \frac{3(1-\mu)yz}{r_1^5} + \frac{3\mu yz}{r_2^5} \\
    H_{1,2} &= -H_{2,1} = 2
\end{align}

Applying the weak form of \ac{PMP}, the optimal thrust direction unit vector is characterized by \cite{Taheri_2017}
\begin{align}
    \hat{\boldsymbol{\alpha}}^* \in \text{arg}\hspace{1mm}\underset{\hat{\boldsymbol{\alpha}}}{\text{min}}\hspace{1mm}H
\end{align}
\noindent which, noting that $uT_{max}/m \geq 0$ and $\lambda_v \ne 0$ in general, gives
\begin{align}
    \Hat{\boldsymbol{\alpha}}^* = -\frac{\boldsymbol{\lambda}_v}{\lambda_v}
    \label{eqn:astar}
\end{align}
\noindent where $\lambda_v$ is the Euclidean norm of $\boldsymbol{\lambda}_v$, and $-\boldsymbol{\lambda}_v$ is defined as the primer vector \cite{Lawden_1964}. Substituting this result into Eq. \eqref{eqn:H1} and defining a switching function $S$ as \cite{Zhang_2015}
\begin{align}
    S = -\frac{c}{m}\lambda_v - \lambda_m + 1,
\end{align}
\noindent the Hamiltonian can now be written as
\begin{align}
    H = \boldsymbol{\lambda}_r^T\mathbf{v} + \boldsymbol{\lambda}_v^T\left[\mathbf{g}(\mathbf{r}) + \mathbf{h}(\mathbf{v})\right] + \frac{uT_{max}}{c}\left(S - \epsilon + \epsilon u\right)
\end{align}

Again applying the weak form of \ac{PMP}, the optimal thrust throttling factor is given by\cite{Zhang_2015}
\begin{equation}
    u^* = \left\{
    \begin{array}{ccc}
        0 & S > \epsilon \\
        \frac{\epsilon - S}{2\epsilon} & -\epsilon \leq S \leq \epsilon \\
        1 & S < -\epsilon
    \end{array}
    \right.
    \label{eqn:ustar}
\end{equation}
Note that the relations between the co-states and optimal control variables $\hat{\boldsymbol{\alpha}}^*$ and $u^*$, given in Eqs. \eqref{eqn:astar} and \eqref{eqn:ustar}, are equivalent to Lawden's primer vector control law when $\epsilon=0$.\footnote{Lawden's primer vector control law is a well-known relationship between the co-states, optimal throttling factor and thrust direction when the indirect optimal control approach is employed to solve minimum-fuel problems \cite{Lawden_1964}.}

To complete the formulation of the trajectory optimization problem, terminal constraints placed on the physical state variables are given by 
\begin{equation}
    \begin{array}{ccc}
        \mathbf{r}(t_i) - \mathbf{r}_i = 0 & \mathbf{v}(t_i) - \mathbf{v}_i = 0 & m(t_i) - 1 = 0 \\
        \mathbf{r}(t_f) - \mathbf{r}_f = 0 & \mathbf{v}(t_f) - \mathbf{v}_f = 0 
    \end{array}
    \label{eqn:state_cons}
\end{equation}
\noindent where $\mathbf{r}_i$, $\mathbf{v}_i$, $\mathbf{r}_f$, and $\mathbf{v}_f$ are the desired initial and final position and velocity vectors, and the initial mass has been scaled to one. Noting that the final mass of the spacecraft is unconstrained, a transversality condition provides the last constraint, which is given by \cite{Bryson_1975}
\begin{align}
    \lambda_m(t_f) = 0
    \label{eqn:costate_con}
\end{align}
\noindent Substituting Eq. \eqref{eqn:astar} into Eq. \eqref{eqn:state_dyn} and combining with Eq. \eqref{eqn:costate_dyn}, a 14-dimensional system of first order ordinary differential equations is obtained
\begin{align}
    \Dot{\mathbf{y}} = \begin{bmatrix} \Dot{\mathbf{x}} \\ \Dot{\boldsymbol{\lambda}} \end{bmatrix} = \begin{bmatrix} 
    \mathbf{v} \\ \mathbf{g}(\mathbf{r}) + \mathbf{h}(\mathbf{v}) - \boldsymbol{\lambda}_vuT_{max}/(\lambda_vm) \\ -uT_{max}/c \\ -\mathbf{G}^T\boldsymbol{\lambda}_v \\ -\boldsymbol{\lambda}_r - \mathbf{H}^T\boldsymbol{\lambda}_v \\ -\lambda_vuT_{max}/m^2
    \end{bmatrix}
    \label{eqn:full_ode}
\end{align}
which, when paired with the optimal throttling factor given in Eq. \eqref{eqn:ustar} and the constraints given in Eqs. \eqref{eqn:state_cons} and \eqref{eqn:costate_con}, fully define the \ac{TPBVP}. The solution of the \ac{TPBVP} is the set of initial co-state variables which satisfy the terminal constraints. Here, it is important to note that solutions of the derived \ac{TPBVP} only imply satisfaction of the first-order necessary conditions of optimality \cite{Russell_2007}. Therefore, without further investigation of the second-order sufficient conditions of optimality, solutions are only guaranteed to be stationary or extremal solutions of the optimal control cost function given in Eq. \eqref{eqn:homotopycost}, also known as candidates of optimality \cite{Prussing_2005}. Furthermore, one should note that multiple local extrema may exist, resulting in multiple solutions of the \ac{TPBVP}.

\section{TPBVP Solution Methodology}
\label{sec:TPBVP}
Solving the derived \ac{TPBVP} requires a good guess of the initial co-state variables $\boldsymbol{\lambda}(t_i)$. The methodology proposed here involves the use of the \ac{PSO} algorithm to first initialize the co-state variables for the minimum-energy problem, corresponding to the perturbation parameter $\epsilon = 1$ in Eq. \eqref{eqn:ustar}. The \ac{PSO}-initialized co-states are then used to seed a single shooting procedure with homotopy continuation, where a trust-region nonlinear solver \cite{Patrick_2020} iteratively discovers the solution to the \ac{TPBVP} for decreasing values of the perturbation parameter until a solution corresponding to $\epsilon = 0$ is found. 

In this Section, the methodology employed to integrate the state, co-state, and \ac{STM} differential equations while performing switching detection is discussed. The proposed method for initializing the co-state variables with \ac{PSO} is also presented, along with an approach to perform single shooting with a homotopy continuation. Along the way, important considerations in regard to programming language selection and computational efficiency are made.

\subsection{Numerical Integration with Switching Detection}
When integrating the state and co-state differential equations, in addition to the \ac{STM} differential equations when single shooting is performed, it is important to ensure that a sufficiently high degree of accuracy is maintained, so that convergence can be achieved. Due to the piecewise continuous function which defines the optimal throttling factor (see Eq. \eqref{eqn:ustar}), integration error can grow about the points at which switching of the throttling factor occurs, and can severely restrict the convergence of a shooting algorithm. However, this can be avoided if the points at which the switching occurs are explicitly determined and numerically integrated up to. Additionally, both the \ac{PSO} co-state initialization and single shooting phases require numerous repeated numerical integrations of the differential equations. Performing numerical integration in a programming language which exhibits high computational performance with a differential equation solving suite that provides high-order Runge-Kutta algorithms and robust event detection is thus desirable. For these reasons, the Julia \cite{Julia} programming language and the Julia-implemented differential equation solving package, \textit{DifferentialEquations.jl} \cite{Rackauckas_2017}, were chosen for all numerical computation and integration in this work.

Julia is a dynamically typed language with syntax similar to that of MATLAB or Python, thereby lending itself to quick development. Historically, a known characteristic of this category of languages is the slow execution time due to the required use of an interpreter, an issue which Julia circumvents through the use of a just-in-time (JIT) LLVM-based compiler. Therefore, although Julia is a dynamic language with easy-to-read syntax, it allows producing software with performance that rivals statically typed languages like C and FORTRAN. Furthermore, \textit{DifferentialEquations.jl} is highly optimized and provides a plethora of modern differential equation solvers and features, including robust event detection schemes, termed \textit{callbacks}, which can employ high-order interpolants to check whether an event has occurred at multiple points backwards in time. Through the application of \textit{DifferentialEquations.jl} and its \textit{continuous callback} feature (i.e. continuous event detection), numerical integration can be performed efficiently by exploiting an adaptive time-stepping method while robustly detecting throttling factor switching to within 64-bit floating point number precision. 

\subsection{Co-State Initialization with PSO}
Beginning with the inception of the \ac{PSO} algorithm detailed in the work of Kennedy and Eberhart \cite{Kennedy_1995}, many have proposed modifications with improved convergence characteristics \cite{Mezura_2011,Pedersen_2010}. The term \ac{PSO} has thus evolved to encompass a class of heuristic global optimization algorithms rather than any one strictly defined algorithm. For the sake of reproducibility, in this work, an implementation of \ac{PSO} in the Julia programming language is employed, which is based on the MATLAB \ac{PSO} algorithm \cite{MATLAB_PSO}.

To employ \ac{PSO} to initialize the co-state variables for the minimum-energy problem ($\epsilon = 1$), a nonlinear, unconstrained optimization problem is formulated as
\begin{align}
    \underset{\boldsymbol{\lambda}(t_i)}{\text{Minimize}} \hspace{3mm} J_{PSO}(\boldsymbol{\lambda}(t_i)) = \mathbf{e}^T(t_f)\mathbf{W}\mathbf{e}(t_f)
    \label{eqn:PSO_Opt_Prob}
\end{align}
\noindent where $\mathbf{W} \in \mathbb{R}^{7\times7}$ is a diagonal weighting matrix and $\mathbf{e}(t_f)$ is the final time boundary condition residual given by
\begin{align}
    \mathbf{e}(t_f) = \begin{bmatrix} \mathbf{r}(t_f) - \mathbf{r}_f \\ \mathbf{v}(t_f) - \mathbf{v}_f \\ \lambda_m(t_f) \end{bmatrix}
\end{align}
\noindent which is computed by integrating the state and co-state differential equations from $t_i$ to $t_f$ as described before, starting from an initial state vector constructed as a concatenation of the constrained initial state and the current guess for the co-state values at the initial time
\begin{align}
    \mathbf{y}(t_i) = \begin{bmatrix} \mathbf{r}_i^T & \mathbf{v}_i^T & 1.0 & \boldsymbol{\lambda}^T(t_i) \end{bmatrix}^T
\end{align}

\noindent Through computational investigation, the implemented single shouting routine, discussed in the following Section, was found to be most sensitive to large residuals in the final position. Therefore, the weighting matrix was defined as
\begin{align}
    \mathbf{W} = \text{Diag}\left(\begin{bmatrix}10 & 10 & 10 & 1 & 1 & 1 & 1 \end{bmatrix}\right)
\end{align}
such that the squared position residuals are weighted by a factor of ten, while the remaining weights are left at unity. 

Due to the stochastic nature of the \ac{PSO} algorithm, some particles may travel to positions that result in trajectories that pass below the surface of a primary body, which is physically non-realizable. Therefore, the integration routine was set up to halt the integration of a given trajectory if the distance between the spacecraft and the surface of a primary body was smaller than zero at any time.
 
When employing \ac{PSO} for co-state initialization, especially when a large swarm size is required and/or the cost function is expensive to evaluate, it is important to take advantage of the parallelizable nature of the algorithm to reduce compute time. In fact, for each iteration of \ac{PSO}, the cost function must be re-evaluated for each particle. In general, the evaluation of the cost function for any given particle is independent of all other particles in the swarm, which makes the process well suited for shared or distributed memory parallelism, which was taken advantage of in this work.

\subsection{Single Shooting with Continuation}
After co-state initialization with \ac{PSO}, a single shooting routine based on a trust-region nonlinear solver provided by \textit{NLsolve.jl} \cite{Patrick_2020} is employed for solving the \ac{TPBVP}. To do so, the nonlinear 7-dimensional vector function which is desired to be solved is defined as
\begin{align}
    \mathbf{Z}(\boldsymbol{\lambda}(t_i),\epsilon) = \mathbf{e}(t_f)
    \label{eqn:shooting_func}
\end{align}
\noindent where a solution corresponds to $\mathbf{e}(t_f) = \mathbf{0}_{7\times 1}$. In other words, it is desired to find the initial co-state vector $\boldsymbol{\lambda}(t_i)$ which, for a given value of $\epsilon \in [0, 1]$ -- integrating Eq. \eqref{eqn:full_ode} from $t_i$ to $t_f$ as described before -- results in a final time boundary condition residual of $\mathbf{e}(t_f) = \mathbf{0}_{7\times1}$. Henceforth, $\mathbf{Z}(\boldsymbol{\lambda}(t_i),\epsilon)$ in Eq. \eqref{eqn:shooting_func} is referred to as the \textit{shooting function}. Seeded by the \ac{PSO}-initialized co-states, the trust-region algorithm is used to iteratively solve Eq. \eqref{eqn:shooting_func} for values of the perturbation parameter decreasing from one to zero according to the continuation law \cite{Zhang_2015}
\begin{align}
    \epsilon_j = \frac{j^2 - 1}{N^2 - 1} \hspace{5mm} j = N, N-1, \dots, 2, 1
    \label{eqn:contlaw}
\end{align}
where $N$ is the number of iterations that are desired to be performed for transitioning the perturbation parameter from one to zero. This value can be chosen somewhat arbitrarily, but it is important to consider that higher values of $N$ result in smaller changes in the perturbation parameter at each iteration, providing a smoother transition from the minimum-energy to the minimum-fuel cost function, albeit at a higher computational cost. Throughout this work, a value of $N=25$ was chosen to balance numerical stability and computational efficiency.

It is important to note that the trust-region algorithm requires that the Jacobian of the shooting function, i.e., $\partial \mathbf{Z}(\boldsymbol{\lambda}(t_i),\epsilon) / \partial \boldsymbol{\lambda}(t_i)$, also be computed. Multiple methods exist to compute this Jacobian (e.g., automatic differentiation and finite difference), among which the \ac{STM} was selected for this work. The \ac{STM} is defined as
\begin{align}
    \boldsymbol{\Phi}(t_f,t_i) = \frac{\partial \mathbf{y}(t_f)}{\partial \mathbf{y}(t_i)}
\end{align}
\noindent To compute the \ac{STM}, the $14^2$ \ac{STM} differential equations must also be integrated along with the state and co-state differential equations in Eq. \eqref{eqn:full_ode} according to
\begin{align}
    \Dot{\boldsymbol{\Phi}}(t,t_i) = \mathbf{F}\boldsymbol{\Phi}(t,t_i) = \frac{\partial \Dot{\mathbf{y}}}{\partial \mathbf{y}}\Big|_t\boldsymbol{\Phi}(t,t_i), \hspace{5mm} \boldsymbol{\Phi}(t_i,t_i) = \mathbf{I}_{14\times14}
\end{align}
\noindent where $\mathbf{F}$ is the Jacobian of Eq. \eqref{eqn:full_ode} evaluated at time $t$, the analytical expression of which is provided by Zhang et al. \cite{Zhang_2015}. Additionally, when discontinuities are present in the shooting function, as is the case for $\epsilon = 0$ when switching occurs (see Eq. \eqref{eqn:ustar}), the \ac{STM} across the discontinuity must also be computed as \cite{Russell_2007}
\begin{align}
    \boldsymbol{\Psi}(t_n) = \frac{\partial \mathbf{y}(t_n^+)}{\partial \mathbf{y}(t_n^-)} = \mathbf{I}_{14\times14} + \left(\Dot{\mathbf{y}}\big|_{t_n^+} - \Dot{\mathbf{y}}\big|_{t_n^-}\right)\left(\frac{\partial S}{\partial \mathbf{y}} \frac{1}{\Dot{S}}\right)
\end{align}
\noindent where $t_n^-$ and $t_n^+$ represent the time immediately before and after the discontinuity occurred, respectively, and
\begin{align}
    \frac{\partial S}{\partial \mathbf{y}} &= \begin{bmatrix} \mathbf{0}_{1\times6} & \frac{c}{m^2}\lambda_v & \mathbf{0}_{1\times3} & -\frac{c}{m}\frac{\boldsymbol{\lambda}_v^T}{\lambda_v} & -1 \end{bmatrix} \\
    \Dot{S} &= \frac{c}{m}\left(\boldsymbol{\lambda}_r + \mathbf{H}^T\boldsymbol{\lambda}_v\right)^T\frac{\boldsymbol{\lambda}_v}{\lambda_v}
\end{align}
\noindent For a trajectory with $M$ discontinuities between $t_i$ and $t_f$, the state transition matrix is then given by \cite{Zhang_2015}
\begin{align}
    \boldsymbol{\Phi}(t_f,t_i) = \boldsymbol{\Phi}(t_f,t_M^+)\boldsymbol{\Psi}(t_M)\boldsymbol{\Phi}(t_{M}^-,t_{M-1}^+)\boldsymbol{\Psi}(t_{M-1})\dots\boldsymbol{\Phi}(t_2^-,t_1^+)\boldsymbol{\Psi}(t_1)\boldsymbol{\Phi}(t_1^-,t_i)
\end{align}

Once the \ac{STM} has been propagated from $t_i$ to $t_f$ along with the state and co-state variables, the Jacobian of the shooting function is constructed from components of the \ac{STM}, such that
\begin{align}
    \frac{\partial \mathbf{Z}(\boldsymbol{\lambda}(t_i),\epsilon)}{\partial\boldsymbol{\lambda}(t_i)} = 
    \begin{bmatrix}
        \frac{\partial\mathbf{r}(t_f)}{\partial\boldsymbol{\lambda}_r(t_i)} & \frac{\partial\mathbf{r}(t_f)}{\partial\boldsymbol{\lambda}_v(t_i)} & \frac{\partial\mathbf{r}(t_f)}{\partial\lambda_m(t_i)} \\
        \frac{\partial\mathbf{v}(t_f)}{\partial\boldsymbol{\lambda}_r(t_i)} & \frac{\partial\mathbf{v}(t_f)}{\partial\boldsymbol{\lambda}_v(t_i)} & \frac{\partial\mathbf{v}(t_f)}{\partial\lambda_m(t_i)} \\
        \frac{\partial\lambda_m(t_f)}{\partial\boldsymbol{\lambda}_r(t_i)} & \frac{\partial\lambda_m(t_f)}{\partial\boldsymbol{\lambda}_v(t_i)} & \frac{\partial\lambda_m(t_f)}{\partial\lambda_m(t_i)} \\
    \end{bmatrix}
    \label{eqn:shootingFuncJac}
\end{align}

A diagram describing the algorithm from co-state initialization to minimum-fuel solution is provided in the Appendix.

\section{Results}
\label{sec:results}

\begin{table}[]
\centering
\caption{Constant Parameters}
\label{tab:const_params}
\begin{tabular}{lcr}
\toprule
Constant & Value & Units \\ 
\midrule
Mass Parameter $\mu$ & $1.21506038\times10^{-2}$ & - \\
Grav. Constant $g_0$ & 9.81 & m/$\text{s}^2$  \\
Time Unit (TU)       & $3.75162997\times10^{5}$ & s \\
Length Unit (LU)     & $3.84400000\times10^{5}$ & km \\
Velocity Unit (VU)   & $1.02462131$ & km/s \\
Mass Unit (MU)       & $m_i$ & kg \\ 
\bottomrule
\end{tabular}
\end{table}

In this Section, case studies are presented to validate the methodology, analyze the performance of the proposed method for co-state initialization, as well as to assess the convergence of the \ac{PSO}-initialized co-states to solutions of the minimum-fuel \ac{TPBVP} for varying particle swarm sizes. The proposed methodology is applied to determine low-thrust minimum-fuel transfers in the following two scenarios in the Earth-Moon system: 1) from a geostationary transfer orbit (GTO) to an $L_1$ halo orbit and 2) from an $L_2$ halo orbit to a $L_1$ halo orbit. Constant parameters used for both scenarios, including the time, length, speed, and mass units employed in the non-dimensionalization of the \ac{CR3BP} dynamics, are provided in Table \ref{tab:const_params}. Note that the mass unit (which is a unique value for a given scenario) is defined as the initial mass of the spacecraft $m_i$. The chosen scenarios provide problems of moderate difficulty to benchmark our approach while allowing for validation and comparison to previous works \cite{Zhang_2015,Aziz_2019,Mingotti_2007,Armellin_2006,Conway_2010,Rasotto_2014}.

Due to the stochastic nature of \ac{PSO}, 100 co-state initialization trials are performed for swarm sizes ranging from 250 to 5000 particles in order for statistical conclusions to be made about the performance of the algorithm. The position of the $k$-th particle in the swarm is defined in the 7-dimensional co-state space as 
\begin{align}
    \boldsymbol{\lambda}_{k}(t_i) = \begin{bmatrix} \boldsymbol{\lambda}_{r,k}^T(t_i) & \boldsymbol{\lambda}_{v,k}^T(t_i) & \lambda_{m,k}(t_i) \end{bmatrix}^T
\end{align}
\noindent The swarm is initialized with particle positions uniformly distributed between upper and lower bounds such that $\boldsymbol{\lambda}_r(t_i) \in [-40,40]$, $\boldsymbol{\lambda}_v(t_i) \in [-2, 2]$, and $\lambda_m(t_i) \in [0, 2]$. It is important to note that these are not constraints placed on the positions of particles, but are rather bounds used to define the distribution which starting particle positions are sampled from. Once the swarm is initialized, particles are then free to travel anywhere in the search space, which is restricted to $\boldsymbol{\lambda}_r(t_i) \in [-100,100]$, $\boldsymbol{\lambda}_v(t_i) \in [-10, 10]$, and $\lambda_m(t_i) \in [0, 10]$. Additionally, convergence criteria are chosen such that a \ac{PSO} trial is halted if the objective function value does not improve by more than $1\times10^{-6}$ over the course of 50 iterations, or if a maximum run time of 30 minutes is met. Throughout the following case studies, the effect of different \ac{PSO} swarm sizes only is analyzed for the discussion to remain tractable. The remaining \ac{PSO} parameters are held at their default values (as defined by MathWorks \cite{MATLAB_PSO}), with the exception of the minimum adaptive neighborhood size fraction, which is set to 0.05.

For all numerical integration performed throughout the case studies, Verner's ``most efficient'' Runge-Kutta 9(8) method was employed, with an absolute and relative tolerance of $10^{-14}$, along with its corresponding $9$-th order interpolant for switching detection \cite{Verner_2010}. All computations were performed on a workstation PC with specifications relevant to this work provided in Table \ref{tab:pc_specs}. At each iteration of \ac{PSO}, the evaluations of the cost function for each particle were computed in parallel using shared memory parallelism on 32 of the available 64 physical cores within the CPU.

\begin{table}[]
\centering
\caption{PC Specifications}
\label{tab:pc_specs}
\begin{tabular}{lr}
\toprule
Specification        & Value  \\ 
\midrule
CPU & AMD Ryzen Threadripper Pro 3995WX \\
RAM & 32 GB DDR4 1600 MHz       \\
Operating System & Windows 10 Enterprise    \\
Julia & Version 1.7.0           \\
\bottomrule
\end{tabular}
\end{table}

\subsection{Scenario 1: GTO to L1 Halo Orbit}

For the GTO to $L_1$ halo orbit transfer problem, the spacecraft was chosen to have an initial mass $m_i = 1500$ kg, a maximum thrust $T_{max} = 10$ N, and a specific impulse $I_{sp} = 3000$ s, as in reference \cite{Zhang_2015}. In the minimum-fuel problem, the time of flight (TOF) was fixed at 8.6404 days. The terminal physical boundary conditions for the \ac{TPBVP}, which identify the initial position of the spacecraft and the final point targeted on the halo orbit, were defined as \cite{Zhang_2015}
\begin{align*}
    \mathbf{r}_i &= \begin{bmatrix} -0.0194885115 &
                                 -0.0160334798 &
                                  0.0 \end{bmatrix}^T\text{ LU} \\
    \mathbf{v}_i &= \begin{bmatrix}  8.9188819237 &
                                 -4.0817936888 & 
                                  0.0 \end{bmatrix}^T\text{ VU} \\
    \mathbf{r}_f &= \begin{bmatrix}  0.8233851820  &
                                  0.0 &
                                 -0.0222775563 \end{bmatrix}^T\text{ LU} \\
    \mathbf{v}_f &= \begin{bmatrix}  0.0 &
                                  0.1341841703 &
                                  0.0 \end{bmatrix}^T\text{ VU}
\end{align*}
\noindent such that the GTO had periapsis and apoapsis altitudes of $h_p=400$ km and $h_a=35864$ km, respectively, and the $L_1$ halo orbit had an out-of-plane amplitude of 8000 km.

\begin{figure}
    \centering
    \begin{subfigure}[b]{0.45\textwidth}
        \centering 
        \includegraphics[width=\textwidth]{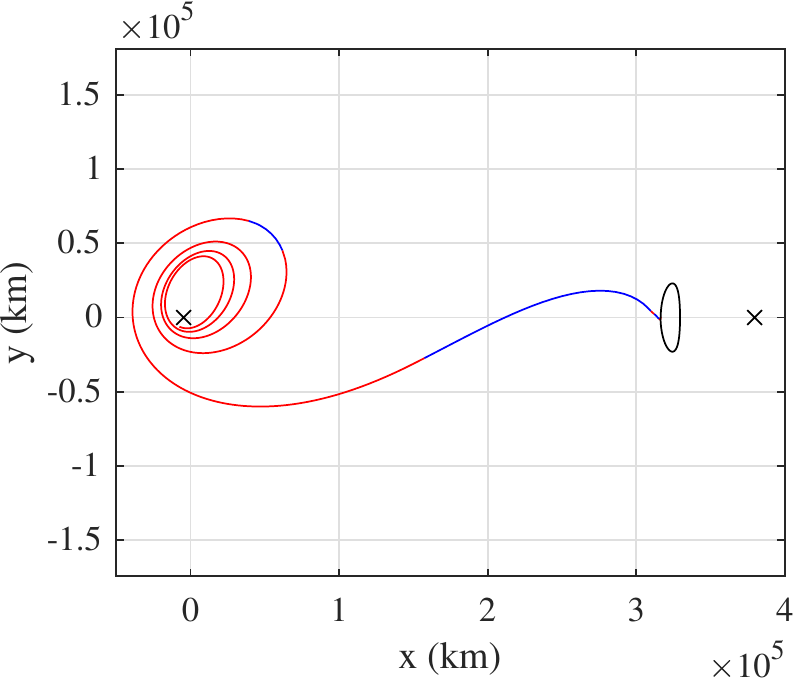}
        \caption{Trajectory A}
        \label{fig:MFTrajA}
    \end{subfigure}
    \hfill 
    \begin{subfigure}[b]{0.45\textwidth}
        \centering 
        \includegraphics[width=\textwidth]{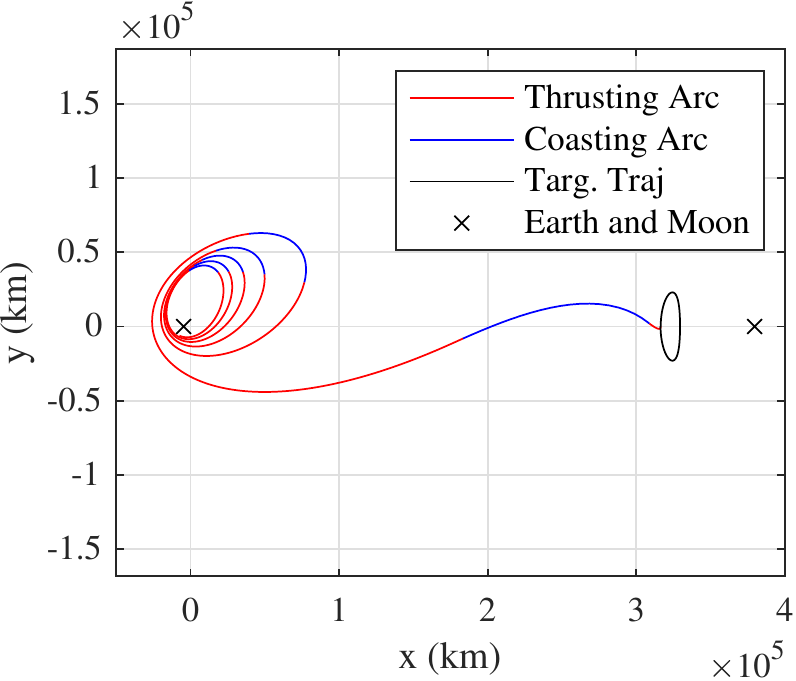}
        \caption{Trajectory B}
        \label{fig:MFTrajB}
    \end{subfigure}
    \begin{subfigure}[b]{0.45\textwidth}
        \centering 
        \includegraphics[width=\textwidth]{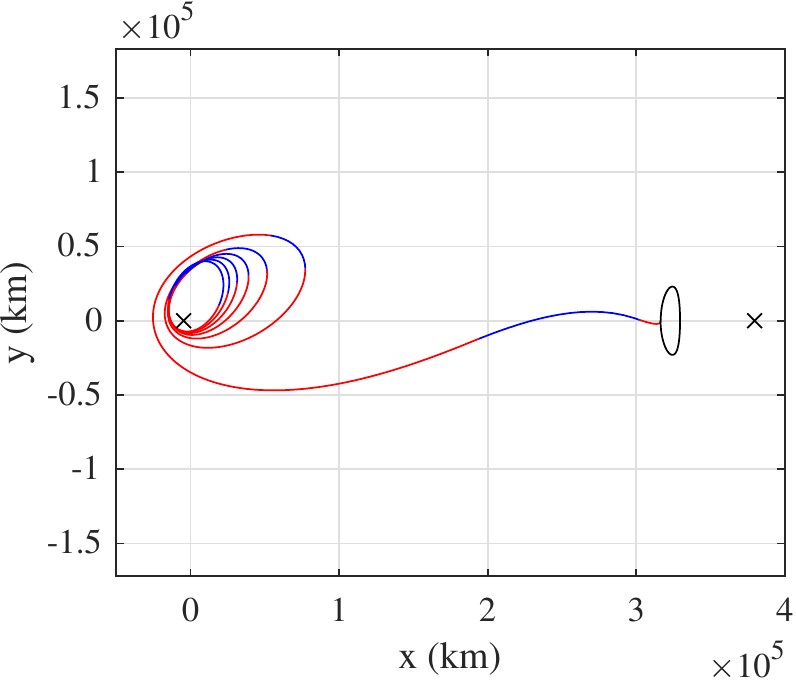}
        \caption{Trajectory C}
        \label{fig:MFTrajC}
    \end{subfigure}
    \hfill 
    \begin{subfigure}[b]{0.45\textwidth}
        \centering 
        \includegraphics[width=\textwidth]{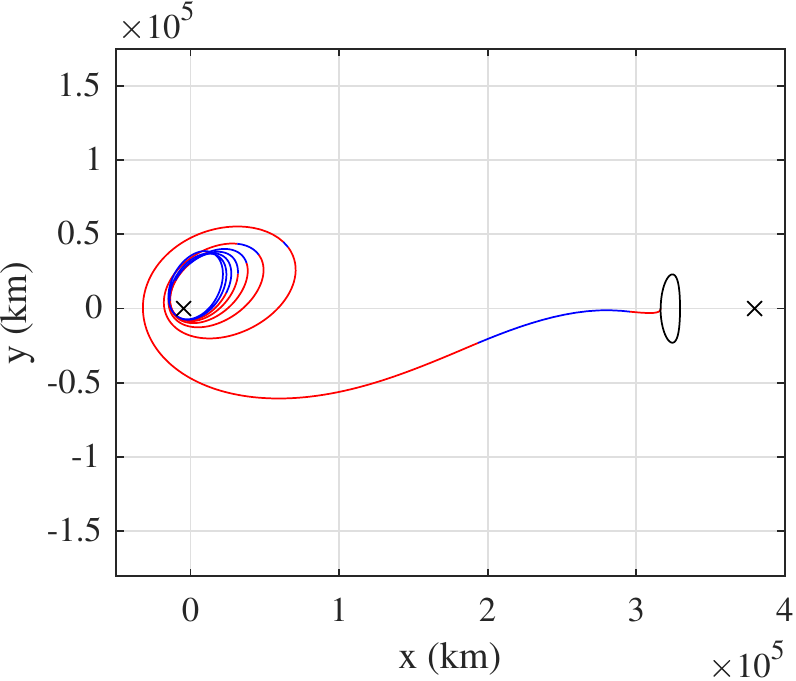}
        \caption{Trajectory D}
        \label{fig:MFTrajD}
    \end{subfigure}
    \begin{subfigure}[b]{0.45\textwidth}
        \centering 
        \includegraphics[width=\textwidth]{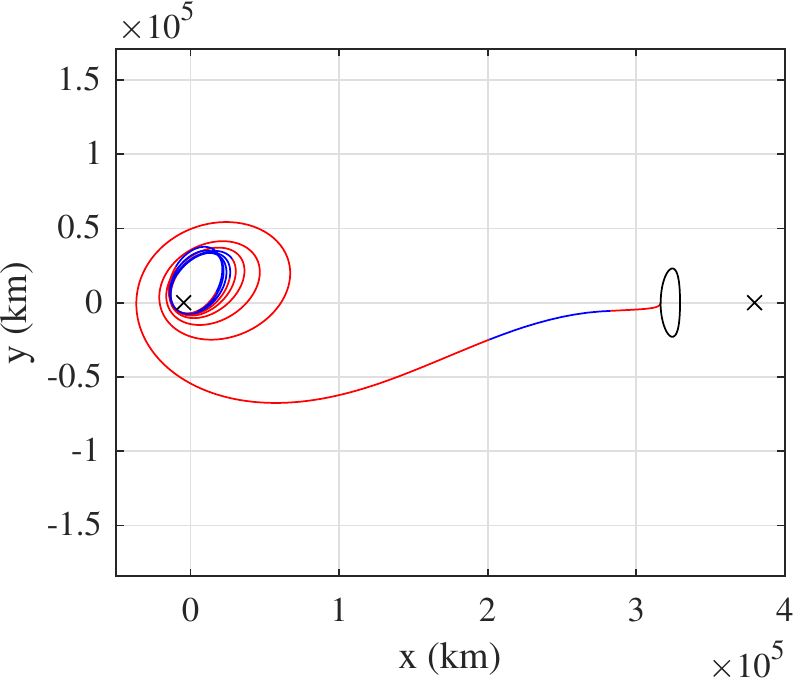}
        \caption{Trajectory E}
        \label{fig:MFTrajE}
    \end{subfigure}
    \hfill 
    \begin{subfigure}[b]{0.45\textwidth}
        \centering 
        \includegraphics[width=\textwidth]{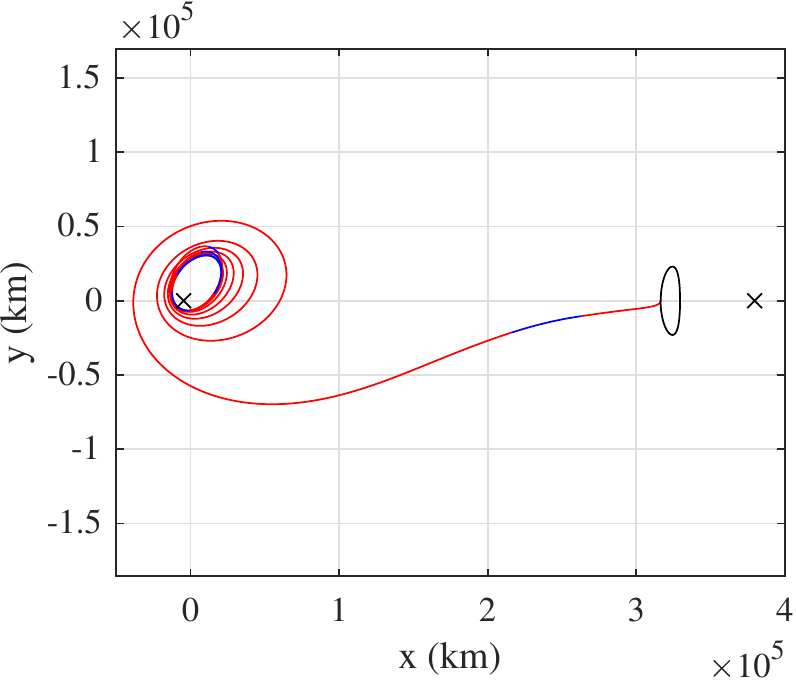}
        \caption{Trajectory F}
        \label{fig:MFTrajF}
    \end{subfigure}
    \caption{Minimum-Fuel Trajectories for Scenario 1}
    \label{fig:MFTrajScenario1}
\end{figure}

\begin{table}[]
\fontsize{7}{7}\selectfont
\centering
\caption{Minimum-Fuel TPBVP Solutions for Scenario 1}
\label{tab:mfsols}
\begin{tabular}{lccccr}
\toprule
Trajectory & $\boldsymbol{\lambda}_r^T(t_i)$ & $\boldsymbol{\lambda}_v^T(t_i)$ & $\lambda_m(t_i)$ & $\Delta m$ (kg) & \# Revs. \\ 
\midrule
A & $\left[23.2524,\hspace{0.1mm}50.6272,-0.08489\right]$ & $\left[-0.1546,\hspace{0.1mm}0.0706,\hspace{0.1mm}-0.0002\right]$ & $0.1385$ & 140.0 & 4 \\
B & $\left[15.6850,\hspace{0.1mm}33.0313,\hspace{0.1mm}-0.0938\right]$ & $\left[-0.1020,\hspace{0.1mm}0.0450,\hspace{0.1mm}-0.0002\right]$ & 0.1334 & 134.4 & 5 \\
C & $\left[7.7397,\hspace{0.1mm}14.5669,\hspace{0.1mm}-0.1239\right]$ & $\left[-0.0466,\hspace{0.1mm}0.0180,\hspace{0.1mm}-0.0001\right]$ & 0.1535 & 139.1 & 6 \\
D & $\left[1.1404,\hspace{0.1mm}-0.6176,\hspace{0.1mm}-0.1535\right]$ & $\left[-0.0006,\hspace{0.1mm}-0.0042,\hspace{0.1mm}-0.0002\right]$ & 0.1935 & 153.5 & 7 \\
E & $\left[-9.4937,\hspace{0.1mm}-25.2095,\hspace{0.1mm}-0.2092\right]$ & $\left[0.0738,\hspace{0.1mm}-0.0401,\hspace{0.1mm}-0.0002\right]$ & 0.2916 & 178.5 & 8 \\
F & $\left[-24.3124,\hspace{0.1mm}-60.4677,\hspace{0.1mm}-0.2883\right]$ & $\left[0.1801,\hspace{0.1mm}-0.0928,\hspace{0.1mm}-0.0002\right]$ & 0.5613 & 219.7 & 9 \\
\bottomrule
\end{tabular}
\end{table}

\subsubsection{Minimum-Fuel TPBVP Solutions}
As discussed before, multiple local extrema of the optimal control cost function, and therefore multiple solutions of the \ac{TPBVP} may exist. A total of six solutions of the minimum-fuel \ac{TPBVP} (corresponding to $\epsilon = 0$) were found. These are labeled Trajectories A through F, and are shown in Figure \ref{fig:MFTrajScenario1}, with thrusting arcs depicted in red and coasting arcs depicted in blue. It is clear that each of the trajectories follows different paths and thrusting and coasting arcs are placed at different locations along the transfers. It is observed that Trajectories A, E, and F (in Figures \ref{fig:MFTrajA}, \ref{fig:MFTrajE}, and \ref{fig:MFTrajF}) exhibit much shorter coasting arcs during either the initial spiraling about the Earth or the final leg of the trajectory when approaching the $L_1$ halo orbit, compared to the remaining trajectories. It is also observed that each trajectory exhibits a unique number of revolutions about the Earth before the final transfer to the halo orbit. 

Table \ref{tab:mfsols} displays the initial co-state variables corresponding to each of the solutions, along with the total change in mass, i.e. fuel required, which is computed as $\Delta m = [1 - m(t_f)]\times1500$ kg. As can be seen by observing the $\Delta m$ column of Table \ref{tab:mfsols}, Trajectory B corresponds to the most optimal minimum-fuel trajectory found, while the remaining solutions correspond to other extremal solutions of the minimum-fuel cost function. The total number of revolutions about the Earth is shown in the final column of Table \ref{tab:mfsols}. The unique number of revolutions for each solution suggests that each trajectory may be locally optimal for a given number of revolutions about the Earth. It is noted that the initial co-states corresponding to Trajectory B, along with $\Delta m$, align with a previously found solution by Zhang et al. \cite{Zhang_2015}.

\subsubsection{Co-State Initialization}
In this Section, the performance of the proposed method for converging to solutions of the derived \ac{TPBVP} which satisfy all boundary conditions is analyzed. For the purpose of this discussion, a trial is considered to have converged if the $\infty$-norm of the \ac{TPBVP} final boundary condition residual is reduced to less than $1\times10^{-10}$, i.e., $||\mathbf{e}(t_f)||_\infty = \max\{|e_i(t_f)| : i=1,2,\dots,7\} < 1\times10^{-10}$. 

\begin{table}[]
\centering
\fontsize{6.5}{7}\selectfont
\caption{Scenario 1 Convergence Characteristics for Different Swarm Sizes}
\label{tab:conv_char}
\begin{tabular}{lcccccccc}
\toprule
Convergence Characteristic & \multicolumn{8}{c}{Swarm Size}  \\
\cmidrule(lr){2-9}
 & 250 & 500 & 750 & 1000 & 2000 & 3000 & 4000 & 5000 \\
\midrule
Min. Fuel \% Converged          & 22  & 42  & 41  & 50  & 63  & 58  & 64  & 54  \\
Avg. Final Obj. Func.           & 7.93   & 5.74   & 5.10   & 4.28   & 3.84   & 3.50   & 2.45   & 2.97   \\
Avg. Time to Converge (min)     & 2.33   & 7.75   & 13.48  & 15.66  & 22.88  & 25.91  & 28.04  & 28.84  \\
Avg. Function Evaluations       & 3.02E5 & 9.24E5 & 1.56E6 & 2.00E6 & 2.82E6 & 3.06E6 & 3.71E6 & 3.58E6 \\
\bottomrule
\end{tabular}
\end{table}

For each of the investigated swarm sizes, Table \ref{tab:conv_char} displays the percentage of trials for which the initialized co-states converged to a solution of the minimum-fuel ($\epsilon=0$) \ac{TPBVP}, the average final value of the \ac{PSO} objective function (see Eq. \eqref{eqn:PSO_Opt_Prob}), and the average time and objective function evaluations before \ac{PSO} convergence. The percentage of converged trials was found to be lowest for the smallest number of particles and to increase with swarm size, reaching a maximum of 64\% converged trials with a swarm size of 4000 particles, before reducing slightly for the largest swarm size. A similar trend of improvement was observed in the average value of the final \ac{PSO} objective function and the number of function evaluations. The average time to converge increases with the swarm size, approaching the maximum time allocated to \ac{PSO} (i.e., 30 minutes per trial) for larger swarm sizes.  

These results suggest that for swarm sizes smaller than 1000 particles, early stagnation of objective function improvement (due to the limited number of particles exploring the solutions space) results in premature convergence of the \ac{PSO} algorithm. For larger swarm sizes, ranging from 2000 to 4000 particles, \ac{PSO} does not suffer from premature convergence and the ability of the algorithm to produce guesses for the initial co-states which converge to solutions of the \ac{TPBVP} is improved greatly. These larger swarm sizes are more effective at exploring the solution space and are able to reduce the \ac{PSO} objective function further, before improvement stalls and the algorithm converges. 
It is also important to note that these larger swarm sizes do come with an additional computational cost: the objective function must be evaluated many more times per \ac{PSO} iteration resulting in longer periods spent waiting for a guess to be generated or prematurely halting optimization depending on the maximum time limit imposed. For a swarm size of 5000 particles, we begin to see the effect of this computational cost: the average time to converge is nearly at the maximum time allocated to \ac{PSO} per trial, indicating that in a majority of cases \ac{PSO} did not converge and was instead halted after 30 minutes of run time (i.e., the \ac{PSO} objective function was still improving when the algorithm was halted). It is expected that this trend would continue for even larger swarm sizes unless more time was allocated for \ac{PSO} co-state initialization.

\subsubsection{Illustration and Analysis of the Solution Generation Process}
Within this Section, the process of generating a solution of the minimum-fuel \ac{TPBVP}, from co-state initialization to single shooting with homotopy continuation, is illustrated and analyzed. The discussion is focused on a single trial that converged to the most fuel-optimal solution found throughout the Scenario 1 case study, i.e., Trajectory B.

\begin{figure}[]
    \centering 
    \begin{subfigure}[b]{0.45\textwidth}
        \centering 
        \includegraphics[width=\textwidth]{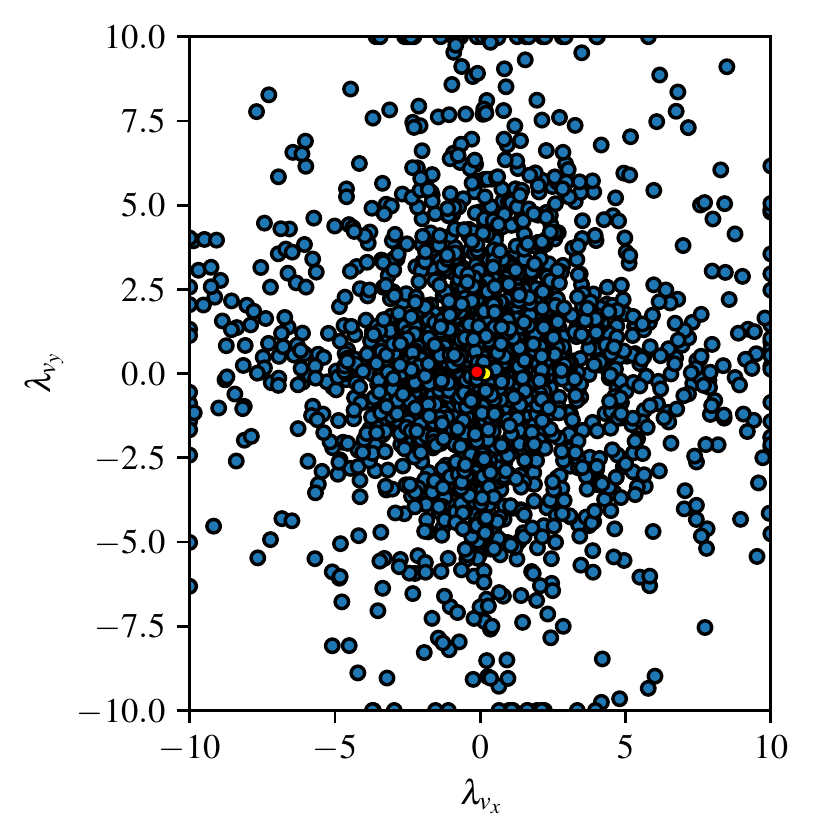}
        \caption{After 5 iterations}
        \label{fig:partInit}
    \end{subfigure} \hspace{10pt}
    \begin{subfigure}[b]{0.45\textwidth}
        \centering 
        \includegraphics[width=\textwidth]{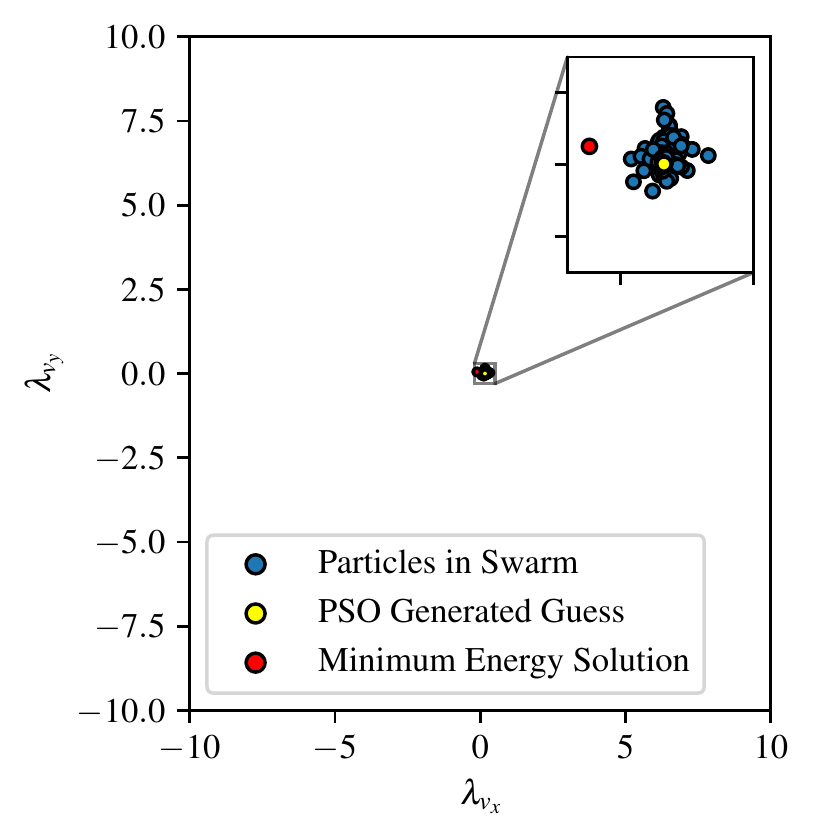}
        \caption{At convergence}
        \label{fig:partFin}
    \end{subfigure}
    \caption{Distribution of \ac{PSO} particles}
    \label{fig:psoDist}
\end{figure}

\begin{figure}
    \centering
    \includegraphics{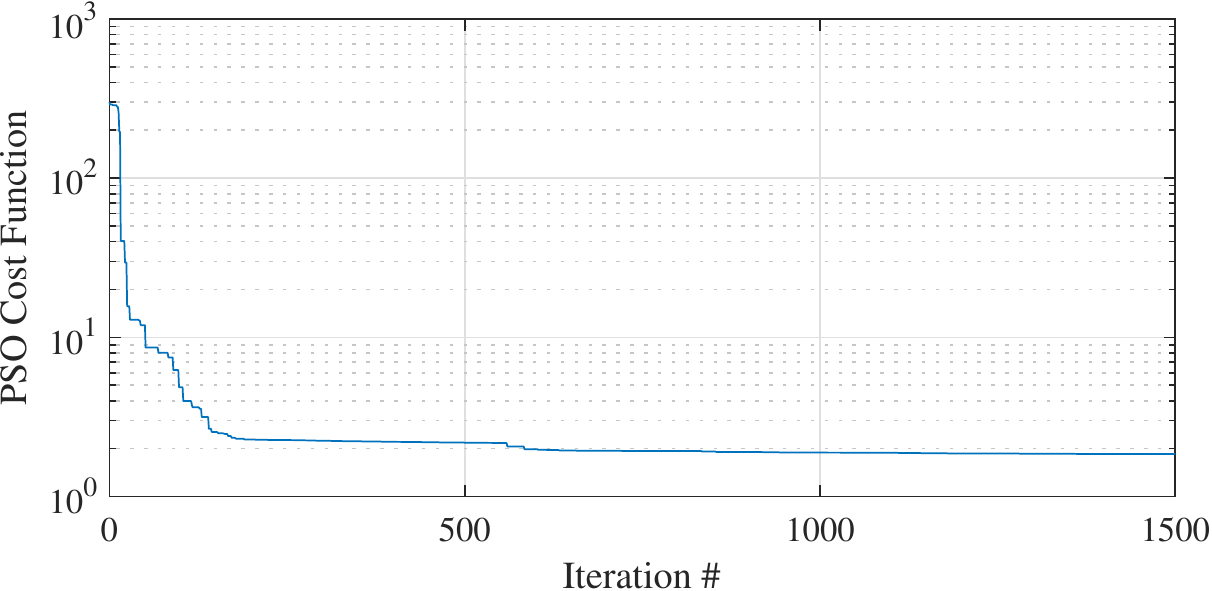}
    \caption{\ac{PSO} objective function versus iteration number}
    \label{fig:psoObj}
\end{figure}

By way of example, 
Figure \ref{fig:psoDist} displays the distribution of 2000 particles (each represented by a blue dot) in the $x$-$y$ velocity co-state space after 5 \ac{PSO} iterations and at \ac{PSO} convergence, along with the location of the \ac{PSO}-generated guess (represented by a yellow dot)  and the minimum-energy solution which the generated guess converged to after performing single shooting (represented by a red dot). Additionally, Figure \ref{fig:psoObj} displays the \ac{PSO} objective function value as the number of \ac{PSO} iterations increases. Observing Figure \ref{fig:partInit}, it can be seen that the particles quickly distributed themselves across the space from their initial placement of $\lambda_{v_x},\lambda_{v_y} \in [-2, 2]$. The particles then spent the next few hundred iterations traveling across the full search space, rapidly decreasing the \ac{PSO} objective function as is shown in Figure \ref{fig:psoObj}. After approximately 200 iterations, once significantly better solutions (i.e., co-states which resulted in a significantly smaller objective function value) were discovered, the rate at which the objective function improved is seen to decrease significantly and the swarm began to contract, spending more time searching the area near the better-known solutions. This slow improvement of the objective function continued until one of the stopping criteria of \ac{PSO} was met. From Figure \ref{fig:partFin}, it can be observed that, upon halting of \ac{PSO}, the particles had traveled to lie within a small region, which happens to be centered about the generated guess. It can also be seen that the generated guess lies near the solution to the minimum-energy problem. This confirms the capability of the algorithm of generating a guess for the co-state variables which lies near a solution of the minimum-energy \ac{TPBVP}.

\begin{figure}[h!]
    \centering
    \includegraphics[width = 0.7\linewidth]{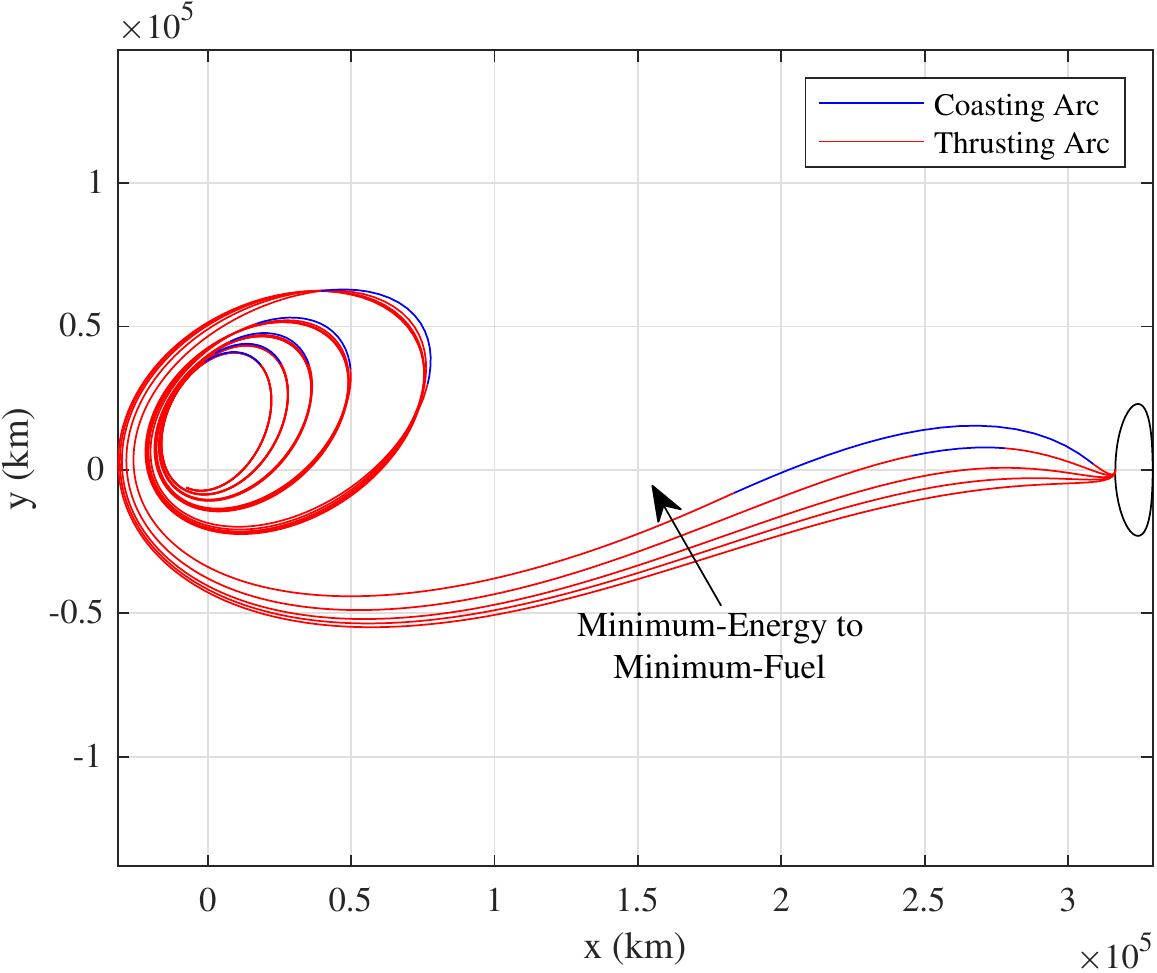}
    \caption{Trajectory B evolution through homotopy continuation}
    \label{fig:trajBevo}
\end{figure}

\begin{figure}[h!]
    \centering
    \includegraphics[width=0.7\linewidth]{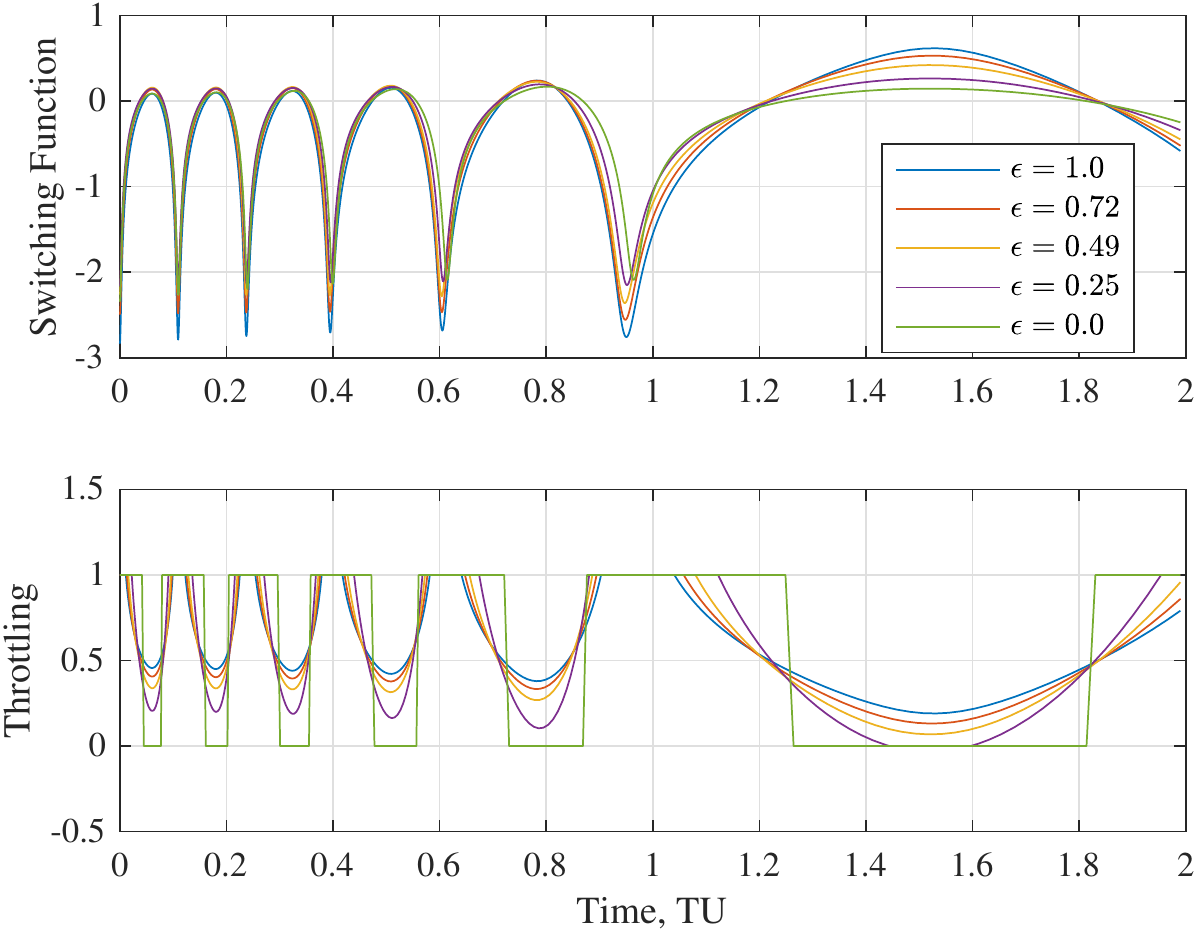}
    \caption{Switching function and throttling factor through Trajectory B continuation}
    \label{fig:switchThrottle}
\end{figure}

After \ac{PSO}-based co-state initialization, the next step in the solution generation process requires employing single shooting to first converge to a solution of the minimum-energy \ac{TPBVP} (i.e., with $\epsilon = 1$) starting from the \ac{PSO}-generated guess, before repeatedly resolving the \ac{TPBVP} while decreasing $\epsilon$ according the Eq. \eqref{eqn:contlaw}. Figure \ref{fig:trajBevo} displays the evolution of Trajectory B through the continuation process as it is transitioned from minimum-energy to minimum-fuel with $\epsilon = 1.0, 0.72, 0.49, 0.25$ and 0.0. One should note that 25 steps in total are taken during the continuation process, and only 5 are selected here for the sake of illustration. Also, Figure \ref{fig:switchThrottle} displays the switching function and throttling factor in time for each of the trajectories displayed in Figure \ref{fig:trajBevo}. As can be seen in Figure \ref{fig:trajBevo}, all trajectories start from and arrive at the same positions as expected. The minimum-energy and minimum-fuel trajectories follow nearly identical paths until the final spiral and transfer to the halo orbit. However, it is observed that coasting arcs appear to increase as the continuation approaches a minimum-fuel trajectory. Figure \ref{fig:switchThrottle} illustrates the gradual introduction of the bang-bang throttling discontinuity throughout the continuation, which is the core strength of the energy-to-fuel homotopy, which widens the convergence radius and improves the probability of \ac{PSO}-generated guesses producing a solution of the minimum-fuel cost function.

\subsection{Scenario 2: L2 Halo Orbit to L1 Halo Orbit}
For the $L_2$ halo orbit to $L_1$ halo orbit transfer problem, the spacecraft was chosen to have an initial mass $m_i=2000$ kg, a maximum thrust $T_{max}=1.5$ N, and a specific impulse $I_{sp}=2000$ s, as in reference \cite{Aziz_2019}. The transfer duration was fixed to 12.7 days and the terminal physical boundary conditions for the \ac{TPBVP}, identifying the initial state of the spacecraft on an $L_2$ halo orbit and the target point on an $L_1$ halo orbit, were defined as \cite{Aziz_2019}
\begin{align*}
    \mathbf{r}_i &= \begin{bmatrix}
        1.1599795702248494 & 0.009720428035815552 & -0.12401864915284157
    \end{bmatrix}^T\text{ LU} \\ 
    \mathbf{v}_i &= \begin{bmatrix}
        0.008477705130550553 & -0.20786307954141953 & -0.010841912833115475
    \end{bmatrix}^T\text{ VU} \\ 
    \mathbf{r}_f &= \begin{bmatrix}
        0.8484736688482315 & 0.00506488863463682 & 0.17343680487577373
    \end{bmatrix}^T\text{ LU} \\ 
    \mathbf{v}_f &= \begin{bmatrix}
        0.005241131023638693 & 0.26343491250951045 & -0.008541420325316247
    \end{bmatrix}^T\text{ VU}
\end{align*} 

\noindent These correspond to initial and target halo orbits with periods of 14.2 and 11.2 days respectively.

\begin{figure}
  \centering
  \begin{subfigure}[b]{0.45\textwidth}
    \centering
    \includegraphics[width=\textwidth]{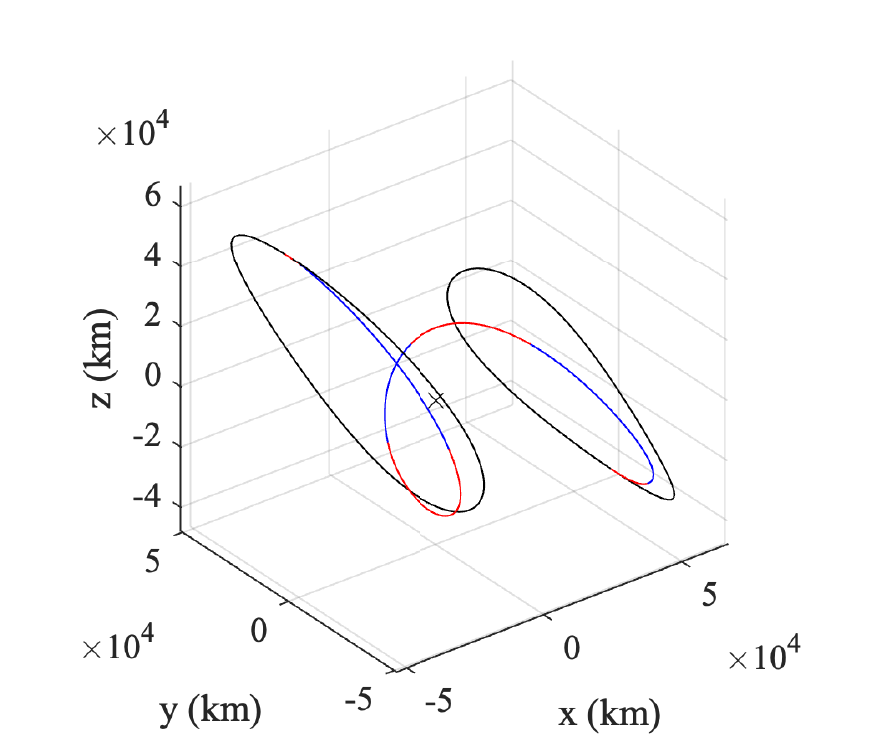}
    \caption{Trajectory $\alpha$ three-dimensional view}
    \label{fig:MFTrajG_3D}
  \end{subfigure} \hspace{10pt}
  \begin{subfigure}[b]{0.45\textwidth}
    \centering
    \includegraphics[width=\textwidth]{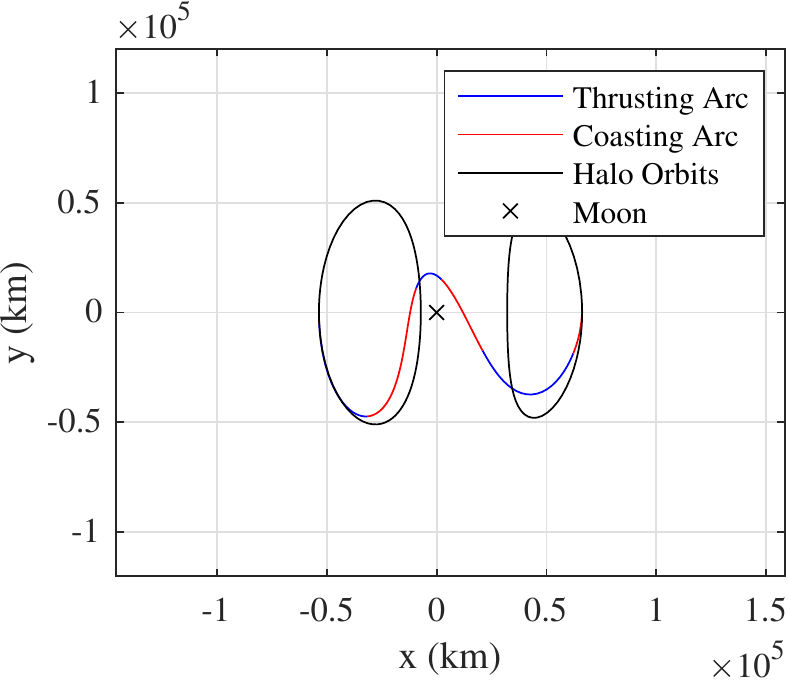}
    \caption{Trajectory $\alpha$ planar projection}
    \label{fig:MFTrajG_XY}
  \end{subfigure}
  \begin{subfigure}[b]{0.45\textwidth}
    \centering
    \includegraphics[width=\textwidth]{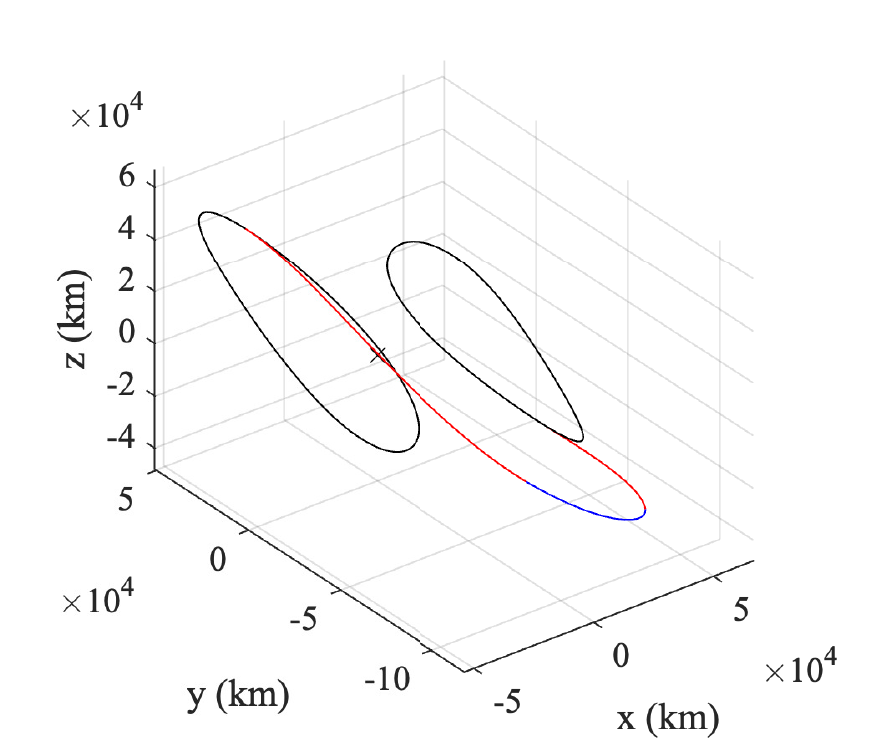}
    \caption{Trajectory $\beta$ three-dimensional view}
    \label{fig:MFTrajH_3D}
  \end{subfigure} \hspace{10pt}
  \begin{subfigure}[b]{0.45\textwidth}
    \centering
    \includegraphics[width=\textwidth]{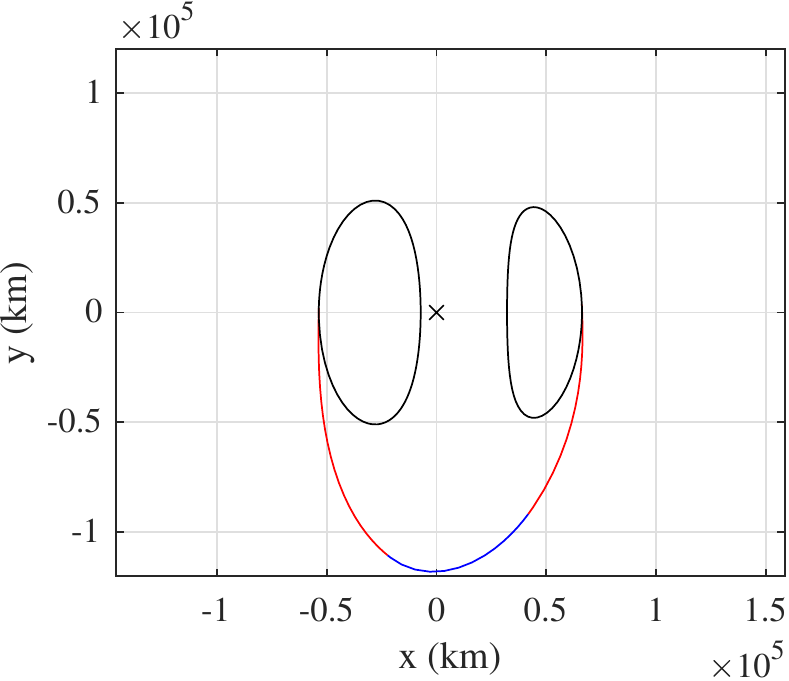}
    \caption{Trajectory $\beta$ planar projection}
    \label{fig:MFTrajH_xy}
  \end{subfigure}
  \begin{subfigure}[b]{0.45\textwidth}
    \centering
    \includegraphics[width=\textwidth]{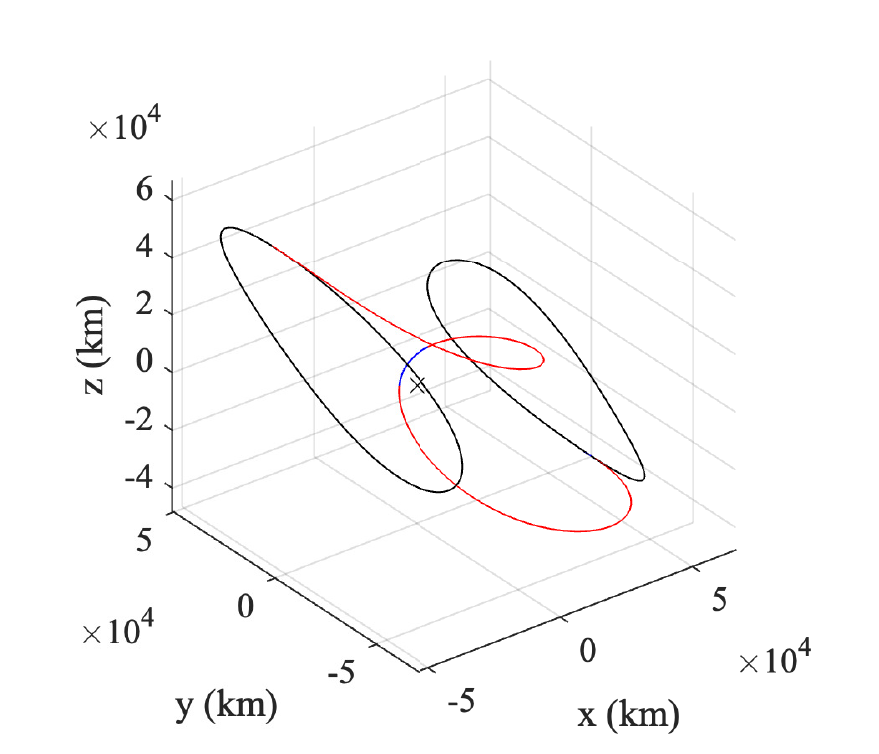}
    \caption{Trajectory $\gamma$ three-dimensional view}
    \label{fig:MFTrajI_3D}
  \end{subfigure} \hspace{10pt}
  \begin{subfigure}[b]{0.45\textwidth}
    \centering
    \includegraphics[width=\textwidth]{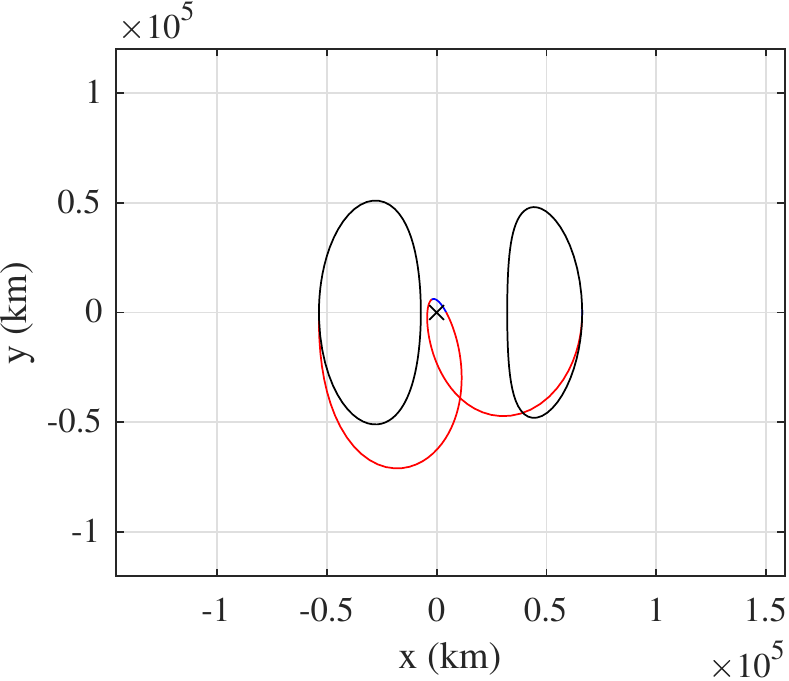}
    \caption{Trajectory $\gamma$ planar projection}
    \label{fig:MFTrajI_xy}
  \end{subfigure}
  \caption{Minimum-Fuel Trajectories for Scenario 2}
  \label{fig:MFTraj_Scenario_2}
\end{figure}

\begin{table}[]
\small
\centering
\fontsize{8}{7}\selectfont
\caption{Minimum-Fuel TPBVP Solutions}
\label{tab:mfsols2}
\begin{tabular}{lcccr}
\toprule
Trajectory & $\boldsymbol{\lambda}_r^T(t_i)$ & $\boldsymbol{\lambda}_v^T(t_i)$ & $\lambda_m(t_i)$ & $\Delta m$ (kg) \\ 
\midrule
$\alpha$ & $\left[0.12603,-0.07665,-0.05635\right]$ & $\left[0.03999, -0.00518, -0.06410\right]$ & $0.02236$ & 35.34 \\
$\beta$ & $\left[-0.01486, 0.01215, -0.07936\right]$ & $\left[0.01015, 0.04457, 0.01256\right]$ & $0.07632$ & 81.28 \\
$\gamma$ & $\left[-0.02195, 0.00659, 0.07490\right]$ & $\left[-0.04314, 0.03615, 0.03842\right]$ & $0.03489$ & 61.27 \\
\bottomrule
\end{tabular}
\end{table}

\subsubsection{Minimum-Fuel TPBVP Solutions}
Similarly to Scenario 1, multiple solutions of the \ac{TPBVP} were discovered by applying the proposed method to solve the $L_2$ halo orbit to $L_1$ halo orbit transfer problem. Figure \ref{fig:MFTraj_Scenario_2} displays the resultant trajectories, both in a three-dimensional view and projected onto the $x$-$y$ plane of the synodic reference frame with the origin shifted to lie at the center of the Moon. It is evident that the three discovered solutions differ drastically. Both Trajectory $\alpha$ and $\beta$ (shown in Figures \ref{fig:MFTrajG_3D}, \ref{fig:MFTrajG_XY}, \ref{fig:MFTrajH_3D}, and \ref{fig:MFTrajH_xy}) appear to take advantage of the lunar gravitational attraction with a close approach of Moon while transferring between halo orbits. Whereas Trajectory $\alpha$ efficiently performs the transfer, visibly requiring the least amount of thrusting among the found solutions, traversing Trajectory $\beta$ requires thrust to be applied for nearly the full duration of the transfer. Employing an entirely different strategy, Trajectory $\gamma$ (shown in Figures \ref{fig:MFTrajI_3D} and \ref{fig:MFTrajI_xy}) does not use a close approach to the lunar surface, and instead thrusts away from the Moon, then follows a coasting arc placed at the greatest distance from the Moon, before thrusting for the final leg of the transfer to insert into the $L_1$ halo orbit.

Table \ref{tab:mfsols2} reports the initial co-state variables corresponding to each of the solutions, along with the fuel required to traverse each trajectory. As expected from the previous discussion, Trajectory $\alpha$ is the most fuel-efficient, requiring just 35.34 kg of fuel, while Trajectory $\gamma$ requires nearly double the fuel mass (i.e., 61.27 kg), and Trajectory $\beta$ is the most fuel expensive (at 81.28 kg). It should be noted that the fuel-optimal trajectory discovered through the application of the proposed method (i.e., Trajectory $\alpha$) agrees with results published by Aziz et al \cite{Aziz_2019}. A slightly lower fuel requirement is found here, which is explained by the fact that the Hybrid Differential Dynamic Programming (HDDP) algorithm applied by Aziz et al. enforced a constant direction of thrust over each integration step and can lag in switching thrust on or off \cite{Aziz_2019}, whereas indirect methods allow thrust to be applied in the optimal direction of the primer vector at every point along the trajectory and enforce throttling on or off at the optimal times.

\subsubsection{Co-State Initialization}

\begin{table}[]
\centering
\tiny
\caption{Convergence Characteristics for Different Swarm Sizes in Scenario 2}
\label{tab:conv_char2}
\begin{tabular}{lcccccccc}
\toprule
Convergence Characteristic & \multicolumn{8}{c}{Swarm Size}  \\
\cmidrule(lr){2-9}
 & 250 & 500 & 750 & 1000 & 2000 & 3000 & 4000 & 5000 \\
\midrule
Min. Fuel \% Converged          & 80   & 86   & 79    & 81   & 67   & 79   & 73   & 75     \\
Avg. Final Obj. Func. ($10^{-3}$)   & 14.6 & 8.89 & 9.57 & 7.37 & 9.33 & 7.56 & 7.61 & 8.30   \\
Avg. Time to Converge (min)     & 2.30    & 3.76    & 4.77    & 6.20    & 10.46   & 11.28   & 12.69   & 14.61     \\
Avg. Function Evaluations       & 7.78E5  & 1.31E6  & 1.58E6  & 2.02E6  & 3.21E6  & 3.25E6  & 3.24E6  & 3.43E6    \\
\bottomrule
\end{tabular}
\end{table}

The performance of the proposed method at converging to solutions of the minimum-fuel $L2$ halo orbit to $L1$ halo orbit transfer \ac{TPBVP} is analyzed in the following. For each of the investigated swarm sizes, Table \ref{tab:conv_char2} displays the percentage of trials that converge to a solution of the minimum-fuel \ac{TPBVP}, along with the final \ac{PSO} objective function value, time to converge, and the number of function evaluations averaged over the 100 trials. 

Overall, the proposed methodology performs more favorably at generating a guess of the initial co-state values for the $L_2$ to $L_1$ transfer problem, compared to Scenario 1. In this case, the best performance was observed when a swarm size of 500 particles was employed, with 86\% of trials converging to a solution of the \ac{TPBVP}. A swarm size of 250 particles performed nearly as well, with 80\% of trials converging to a \ac{TPBVP} solution while only requiring 2.3 minutes on average to generate the guess. Swarm sizes increasing up to 1000 particles continued to produce quality guesses frequently, converging in 79\% and 81\% of trials for 750 or 1000 particles respectively, albeit at a higher computational cost. Interestingly, for a swarm size greater than 1000 particles (aside from 3000 particles), co-state initialization performance appeared to taper off, whereas the average time to converge in all cases was well below the time limit of 30 minutes, indicating that allocation of time for co-state initialization was not the issue. This contradicts results from Scenario 1, which suggested that a larger swarm size should be preferred as long as enough time is allocated for co-state initialization. Furthermore, on average, larger swarm sizes successfully reduced the \ac{PSO} objective function further compared to the smaller, more often successful swarm sizes. This is counter-intuitive, as a smaller average final objective function value indicates that the generated co-state guesses are closer to satisfying the final time boundary condition constraints and are therefore expected to converge to a solution of the \ac{TPBVP} with a higher frequency, as was observed in Scenario 1.

It was discovered that multiple local minima of the \ac{PSO} cost function exist for the problem investigated in Scenario 2. The necessary condition for a minimum of the \ac{PSO} cost function is given by
\begin{equation}
    \left(\frac{\partial J_{PSO}}{\partial \boldsymbol{\lambda}(t_i)}\right)^T = 2\left(\frac{\partial \mathbf{e}(t_f)}{\partial \boldsymbol{\lambda}(t_i)}\right)^T\mathbf{W}\mathbf{e}(t_f) = \mathbf{0}_{7\times1} 
\end{equation}
which is satisfied either 1) when the final boundary condition residual $\mathbf{e}(t_f)$ is reduced to the zero vector or 2) when $\left(\partial \mathbf{e}(t_f)/\partial \boldsymbol{\lambda}(t_i)\right)^T$ becomes rank deficient and $\mathbf{W}\mathbf{e}(t_f)$ lies in its null space. Clearly, global minima of the \ac{PSO} cost function correspond to a final boundary condition residual of zero, and local minima of the \ac{PSO} cost function exist in regions where the Jacobian of the final boundary condition residual with respect to the initial co-state vector is rank deficient.

It is expected that the counter-intuitive phenomenon of larger \ac{PSO} swarm sizes resulting in lower convergence rates observed in Scenario 2 is due to \ac{PSO} converging to regions near local minima of the \ac{PSO} cost function more frequently. In trials when \ac{PSO} 
produced a guess which did not converge to a solution through single shooting, the \ac{PSO}-generated guess was confirmed to lie in a region where $\partial \mathbf{e}(t_f)/\partial \boldsymbol{\lambda}(t_i)$, and therefore the Jacobian of the shooting function with respect to the initial co-states (see Eq. \eqref{eqn:shootingFuncJac}), is rank deficient. This rank deficiency then results in updates of Newton's method in the trust-region algorithm which produce little to no improvement during single shooting and therefore fail to converge to a solution of the \ac{TPBVP}.

\begin{figure}[hbt!]
  \centering
  \begin{subfigure}[b]{0.45\textwidth}
    \centering
    \includegraphics[width=\textwidth]{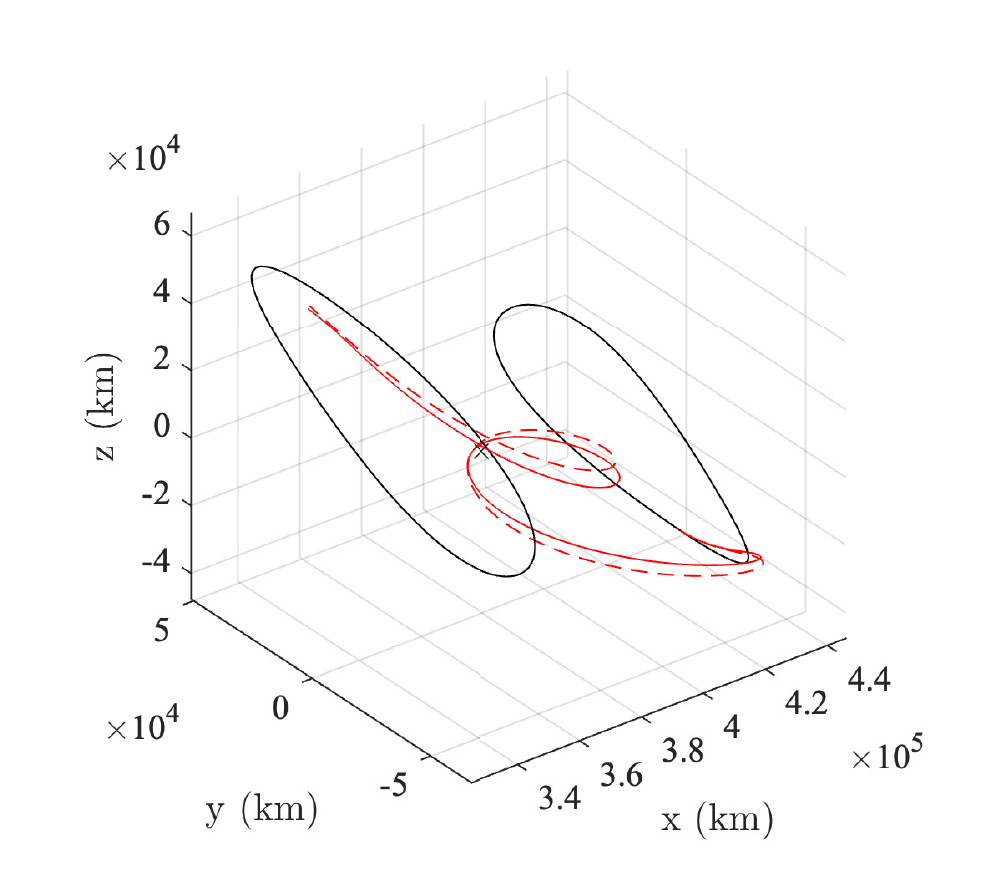}
    \caption{Three-dimensional view}
    \label{fig:Failed_3D}
  \end{subfigure} \hspace{10pt}
  \begin{subfigure}[b]{0.45\textwidth}
    \centering
    \includegraphics[width=\textwidth]{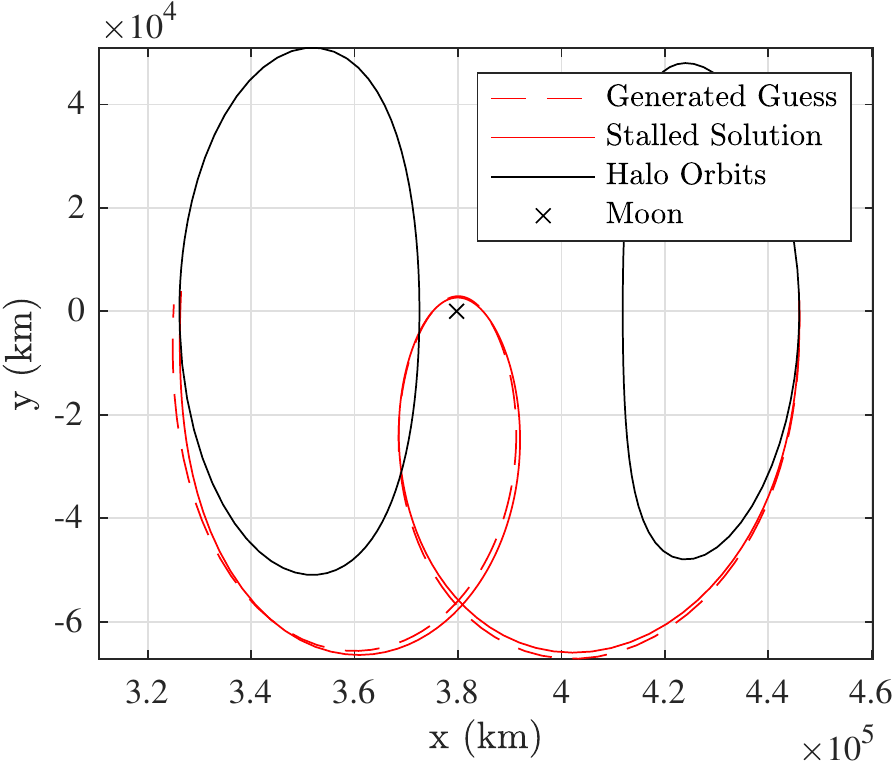}
    \caption{Planar projection}
    \label{fig:Failed_XY}
  \end{subfigure}
  \caption{Initial guess and stalled solution for minimum-energy trajectory}
  \label{fig:Failed}
\end{figure}

An example of a trial that encountered this difficulty is shown in Figure \ref{fig:Failed}. Represented with a dashed line is the trajectory corresponding to the \ac{PSO}-generated guess of

\begin{equation*}
    \boldsymbol{\lambda}_{PSO}(t_i) = \begin{bmatrix}
       -0.75848 & 0.28231 & -9.76567 & -1.27696 & 2.08707 & -3.41462 & 4.57662
    \end{bmatrix}^T
\end{equation*}

\noindent After performing single shooting and failing to converge, the trajectory shown as a solid red line corresponds to the values of the initial co-states at which single shooting stalled, which are given by

\begin{equation*}
    \boldsymbol{\lambda}(t_i) = \begin{bmatrix}
        -0.15992 & -0.23058 & -5.42275 & -0.57888, 0.97406 & -1.98589 & 2.11896
    \end{bmatrix}^T
\end{equation*}

\noindent It is important to note that both trajectories are entirely shown in red because no coasting arcs were exploited by either; however, it should be recalled that the minimum-energy problem does not require optimal trajectories to maintain a bang-bang thrust profile (see Eq. \eqref{eqn:ustar} with $\epsilon=1.0$). The $\infty-$norm of the \ac{TPBVP} final boundary condition residual was $||\mathbf{e}(t_f)||_\infty = 3.95\times10^{-2}$ and $||\mathbf{e}(t_f)||_\infty = 9.72\times10^{-3}$ for the generated guess and the co-states returned after stalled single shooting respectively, which does confirm slight improvement during single shooting. Furthermore, it was found that both the generated guess and single shooting-returned co-states result in a Jacobian of the shooting function with respect to the initial co-states of rank six, which corresponds to rank deficiency. This example illustrates a trial where the final boundary condition residuals were nearly satisfied for both cases but rank deficiency of the Jacobian of the shooting function prevented convergence to a solution of the \ac{TPBVP}.

\section{Conclusion}
This paper proposed a method for initializing the co-state variables for solving low-thrust minimum-fuel trajectory optimization problems with the indirect shooting approach by employing \ac{PSO} paired with an energy-to-fuel homotopy technique. The method was applied successfully to solve two low-thrust transfer problems in the Earth-Moon system. It was demonstrated that the methodology is able to successfully generate guesses for the initial co-state variables which converge to a solution for both investigated scenarios. The resulting minimum-fuel trajectories were validated by comparison with published solutions of the optimal low-thrust transfer for both scenarios. Analysis was performed to determine the effect that varying the number of particles in the swarm has on the performance of the proposed method. It was found that for problems involving several low-thrust spirals about a primary body, a larger swarm size of 2000 to 4000 particles produced the highest quality guesses for the initial co-state variables, whereas, for problems involving short transfers between orbits in cislunar space, a smaller swarm size of 500 particles is preferred and is capable of rapidly generating solutions in less than 4 minutes. It should be noted that, due to the heuristic nature of the \ac{PSO} algorithm, multiple solutions to the derived minimum-fuel \ac{TPBVP} were found, which correspond to extremal solutions of the minimum-fuel objective function, providing interesting physical insight into the investigated trajectory optimization scenarios. Considerations on the existence of local minima and Jacobian rank deficiency were also discussed.

The primary benefits of applying \ac{PSO} for minimum-fuel co-state initialization are the simplicity and versatility of the approach. Once the \ac{TPBVP} is derived through the application of \ac{COV} and \ac{PMP}, little additional work is required to employ \ac{PSO} to initialize the co-state variables for a range of trajectory optimization problems in similar dynamical environments. Additionally, \ac{PSO} removes the need for a guess to be provided by the end-user of the algorithm, an advantageous feature for preliminary mission design efforts, when a large search for many candidate trajectories is required. A large search for candidate trajectories could also be accelerated drastically through the exploitation of a high-performance computing environment, as the methodology developed herein lends itself well to parallel computing and could be distributed across many CPUs to generate solutions for many initial and target states, flight times, and spacecraft parameters. Furthering this cause, the proposed methodology was shown to discover multiple solutions of the derived \ac{TPBVP}, most of which could be considered candidate trajectories during preliminary mission design efforts once local optimality is verified through an investigation of the second-order sufficient conditions of optimality. Finally, while not investigated in this work, the multiple solutions discovered through the application of the proposed methodology could be improved through further homotopy continuation on the fixed flight time or thrust magnitude to produce even more candidate trajectories. 

\section*{Acknowledgments}
This material is based upon work supported by the National Science Foundation Graduate Research Fellowship under Grant No. 88983. Grant Hecht also acknowledges reception of the NASA Space Grant Fellowship. 
        
\bibliographystyle{elsarticle-num} 
\bibliography{references}

\clearpage
\section*{Appendix: Algorithm for PSO Co-State Initialization}
\vspace{12mm}
\begin{figure}[h!]
    \centering
    \includegraphics[width=14cm]{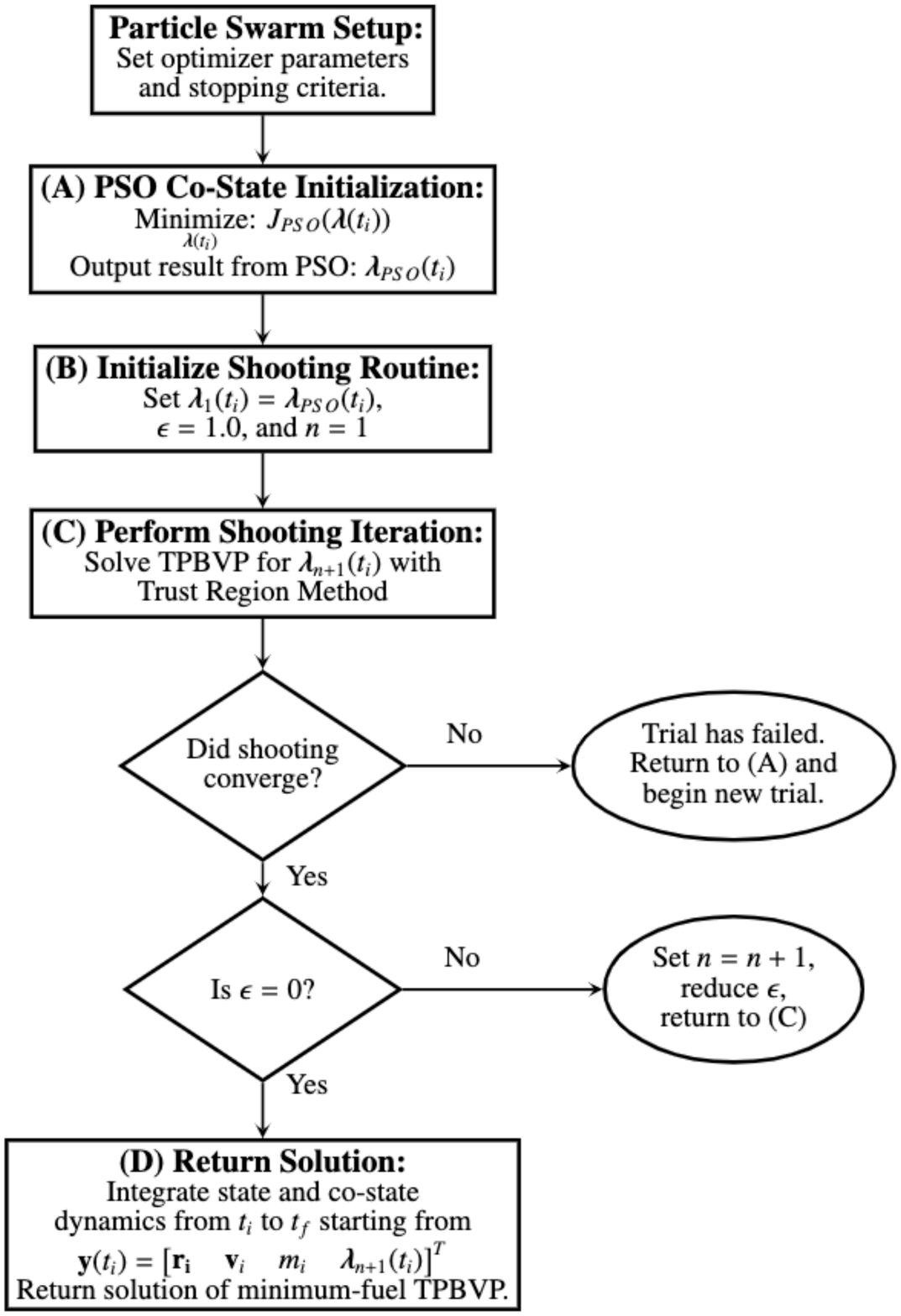}
    \label{dig:AlgDigram}
\end{figure}

\clearpage
\section*{Vitae}
 \begin{figure}[h!]
 \includegraphics[width=3.5cm]{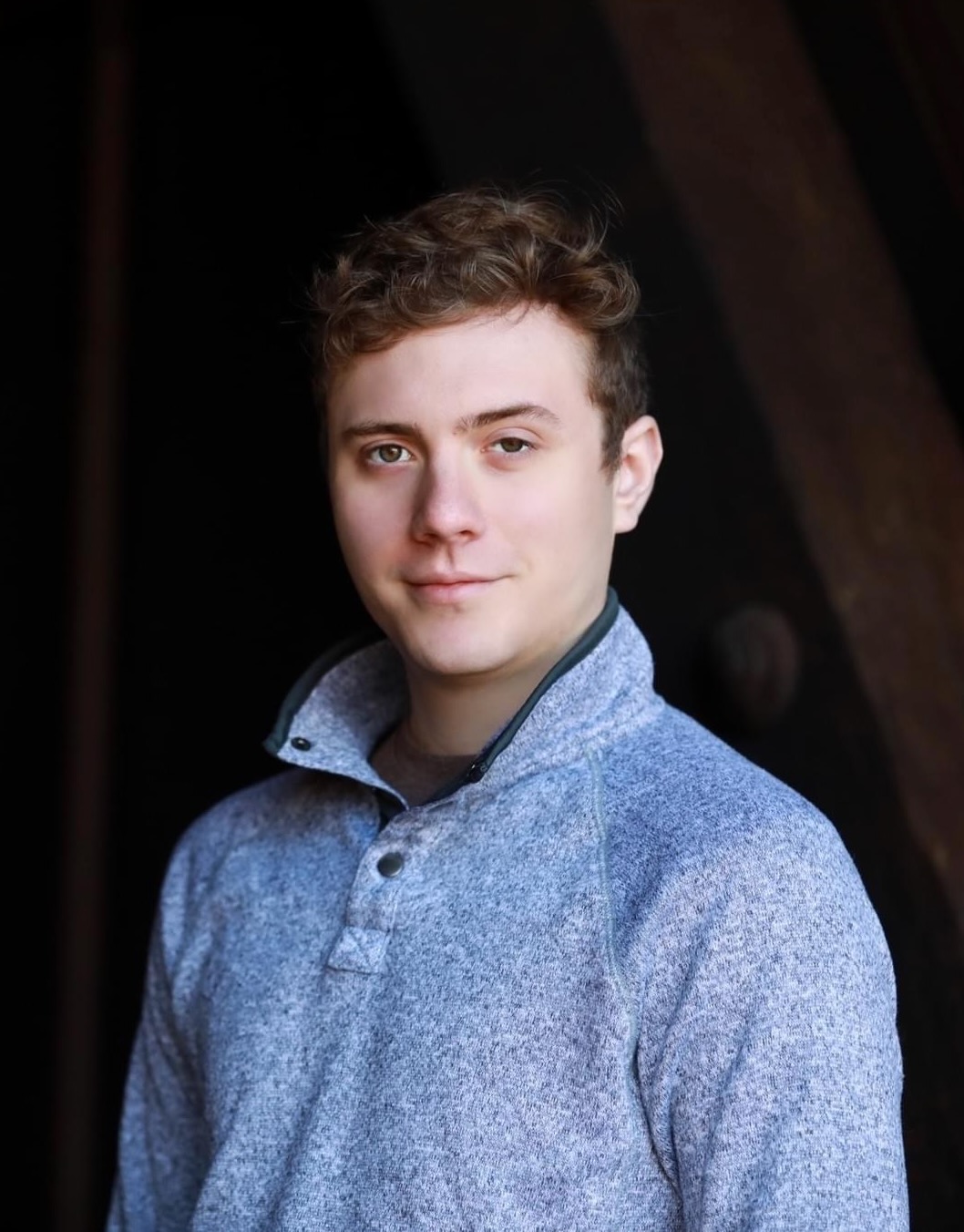}
 \end{figure}
 Grant R. Hecht is a Ph.D. candidate in the Department of Mechanical and Aerospace Engineering at the University at Buffalo. He is also an NSF GRFP fellow and a NASA Pathways Student. He holds a B.Sc. in Aerospace Engineering from the Missouri University of Science and Technology. His main research interests are in spacecraft trajectory optimization, formation guidance, and inertial navigation.

\begin{figure}[h!]
\includegraphics[width=3.5cm]{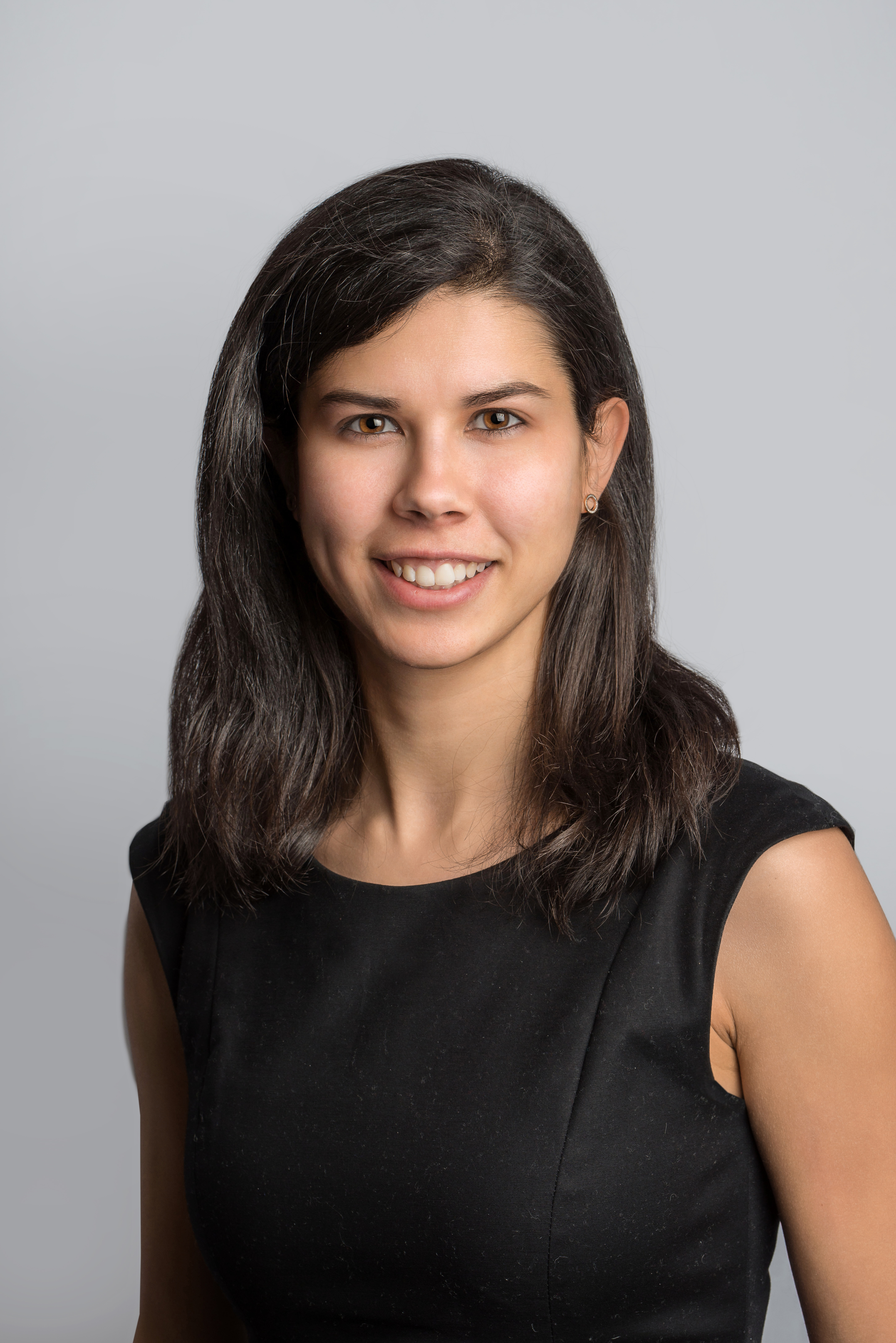}
\end{figure}
Dr. Eleonora M. Botta is an assistant professor in the Department of Mechanical and Aerospace Engineering at the University at Buffalo, NY. She holds a Ph.D. in Mechanical Engineering from McGill University, a B.Sc. in Aerospace Engineering and a M.Sc. in Space Engineering from Politecnico di Milano, and a M.Sc. in Aerospace Engineering from Politecnico di Torino. Her main research interests are the dynamics and control of tethered space systems, spacecraft trajectory optimization, and space situational awareness. She is a member of the AAS Space Flight Mechanics Committee and the vice-chair of the AIAA Space Tethers Technical Committee.

\end{document}